\documentclass[a4paper]{article}

\usepackage[english]{babel}
\usepackage[utf8x]{inputenc}
\usepackage[T1]{fontenc}
\usepackage[a4paper,top=3cm,bottom=2cm,left=3cm,right=3cm,marginparwidth=1.75cm]{geometry}
\usepackage{amsmath}
\usepackage{graphicx}
\usepackage{todonotes}
\usepackage{listings}
\usepackage{color} 
\definecolor{mygreen}{RGB}{28,172,0} 
\definecolor{mylilas}{RGB}{170,55,241}
\usepackage{graphicx}
\usepackage[small]{caption}
\usepackage{subcaption}
\usepackage[makeroom]{cancel}
\usepackage{listings}
\usepackage{amsmath}
\usepackage{amsthm}
\usepackage{amsfonts}
\usepackage{float}
\usepackage{framed}
\usepackage{cite}
\usepackage{algorithm}
\usepackage{pdfpages}
\usepackage{mathtools}
\usepackage{bm}
\usepackage{authblk}
\usepackage{multirow}
\usepackage{authblk}
\usepackage{pbox}

\DeclareMathOperator{\arctantwo}{arctan2}
\newcommand{\norm}[1]{\left\lVert#1\right\rVert}

\newtheorem{theorem}{Theorem}

\newtheorem{remark}{Remark}

\providecommand{\keywords}[1]{\textit{Keywords and phrases: } #1}

\begin{document}
	
\title{A micro-macro Markov chain Monte Carlo method for molecular dynamics using reaction coordinate proposals II: indirect reconstruction}

\author[1]{Hannes Vandecasteele}
\author[1]{Giovanni Samaey}
\affil[1]{KU Leuven, Department of Computer Science, NUMA Section, Celestijnenlaan 200A box 2402, 3001 Leuven, Belgium, hannes.vandecasteele@kuleuven.be}
\date{\today}
	
\maketitle

\begin{abstract}
We introduce a new micro-macro Markov chain Monte Carlo method (mM-MCMC) with indirect reconstruction to sample invariant distributions of molecular dynamics systems that exhibit a time-scale separation between the microscopic (fast) dynamics, and the macroscopic (slow) dynamics of some low-dimensional set of reaction coordinates. The algorithm enhances exploration of the state space in the presence of metastability by allowing larger proposal moves at the macroscopic level, on which a conditional accept-reject procedure is applied. Only when the macroscopic proposal is accepted, the full microscopic state is reconstructed from the newly sampled reaction coordinate value and is subjected to a second accept/reject procedure. The computational gain stems from the fact that most proposals are rejected at the macroscopic level, at low computational cost, while microscopic states, once reconstructed, are almost always accepted. This paper discusses an indirect  method to reconstruct microscopic samples from macroscopic reaction coordinate values, that can also be applied in cases where direct reconstruction is cumbersome.  The indirect reconstruction method generates a microscopic sample by performing a biased microscopic simulation, starting from the previous microscopic sample and driving the microscopic state towards the proposed reaction coordinate value. We show numerically that the mM-MCMC scheme with indirect reconstruction can significantly extend the range of applicability of the mM-MCMC method.
\end{abstract}


\section{Introduction} \label{sec:introduction}
Countless systems in chemistry and physics consist of a large number of microscopic particles, of which all positions are collected in the system state $x\in\mathbb{R}^d$, with $d$ the (high) dimension of the system~\cite{leimkuhler2016molecular}. The dynamics of such systems is usually governed by a potential energy $V(x)$ and Brownian motion $W_t$, for instance through the overdamped Langevin dynamics
\begin{equation}\label{eq:overdamped_langevin}
	dX_t = -\nabla V(X_t)dt+\sqrt{2\beta^{-1}}dW_t,
\end{equation} 
in which $X_t$ represents the time-dependent state of an individual realisation of the dynamics, and $\beta$ is the inverse temperature.  One common computational task for such systems is sampling their time-invariant distribution, i.e., the Gibbs-measure
\begin{equation} \label{eq:mu}
d\mu(x) = Z_V^{-1} \exp\left(-\beta V(x)\right)  dx,
\end{equation}
with $Z_V$ the normalization constant and $dx$ the Lebesgue measure. The standard Metropolis-Hastings algorithm~\cite{metropolis1953equation,hastings1970monte} faces issues when there is a large time-scale separation between the fast dynamics of the full, high-dimensional (microscopic) system and the slow behaviour of some low-dimensional (macroscopic) degrees of freedom. Then, for stability reasons, simulating the microscopic dynamics~\eqref{eq:overdamped_langevin} requires taking time steps on the order of the fastest mode of the system, limiting the size of proposal moves and slowing down exploration of the full state space. In particular, when the potential $V$ contains multiple local minima, standard MCMC methods can remain stuck for a long time in these minima. This phenomenon is called metastability.

There exist several techniques to accelerate sampling in such a context, for instance the parallel replica dynamics~\cite{voter1998parallel,voter2002extending,le2012mathematical,lelievre2016partial}, the adaptive multilevel splitting method~\cite{cerou2007adaptive} and kinetic Monte Carlo~\cite{voter2007introduction}, or (adaptive) biased forcing methods~\cite{wang2001efficient,henin2004overcoming,darve2001calculating, stoltz2010free}.
This paper is the second part of a two-part work on a novel Markov-chain Monte Carlo (MCMC) method in which proposal moves are based on an approximate effective dynamics for some low-dimensional set of reaction coordinates. In the first part of this work \cite{vandecasteele2020direct}, we introduced such a micro-macro Markov chain Monte Carlo method (mM-MCMC) in which a new sample is created following a three-step procedure: (i) \emph{restriction}, i.e., computation of the reaction coordinate value $z$ associated to the current microscopic sample $x$; (ii) a \emph{macroscopic MCMC step}, i.e., sampling a new value $z'$ of the reaction coordinate based on an approximate effective dynamics; and (iii) \emph{reconstruction}, i.e., creation of a microscopic sample $x'$ that is consistent with the sampled macroscopic reaction coordinate value. Both step (ii) and step (iii) contain an accept/reject procedure.  If the proposed reaction coordinate value is accepted during step (ii), we continue with step (iii); otherwise we return to step (i). During step (iii), a second accept/reject step ensures that the microscopic samples indeed sample the Gibbs measure~\eqref{eq:mu}. The algorithm in \cite{vandecasteele2020direct} requires a reconstruction distribution of microscopic samples conditioned upon a given reaction coordinate value.  We call the corresponding method \emph{mM-MCMC with direct reconstruction}.  

The mM-MCMC method with direct reconstruction is unbiased regardless of the choices in the effective dynamics and the reconstruction distribution. However, the computational advantage of the method crucially depends on both choices, for two reasons \cite{vandecasteele2020direct}. First, one needs to ensure that the fastest modes are not present at the reaction coordinate level, such that larger moves are possible than at the microscopic level, enhancing the exploration of the phase space. Second, the scheme should be constructed such that most rejected proposals are already rejected at the reaction coordinate level, i.e., without ever having to perform the costly reconstruction of a corresponding microscopic sample. In particular, the acceptance rate of the reconstructed microscopic samples should be close to $1$.  We have shown that, when these conditions are met, the mM-MCMC method is able to obtain an efficiency gain over the standard Markov chain Monte Carlo method that is proportional to the time-scale separation present in the system \cite{vandecasteele2020direct}. 

In this second paper, we propose an alternative (indirect) reconstruction procedure to perform the reconstruction step after a reaction coordinate has been accepted at the macroscopic level. Often, a direct reconstruction of a microscopic sample for a given value of the reaction coordinate, via sampling from a reconstruction distribution, is computationally costly and potentially cumbersome, see \cite{vandecasteele2020direct}. Indeed, we need to sample from a reconstruction distribution defined on a (possibly highly) non-linear sub-manifold. 
To overcome this issue, we introduce a general \emph{indirect reconstruction} scheme to approximately reconstruct a microscopic sample that lies close to the sub-manifold $\Sigma(z')$ of microscopic samples with reaction coordinate value $z'$.
The main idea of indirect reconstruction is to perform time integration of a strongly biased stochastic process that drives the reaction coordinate value of the last accepted microscopic sample towards the reaction coordinate value $z'$ sampled at the macroscopic level. We then define the reconstructed microscopic sample $x'$ as the microscopic state that is obtained at the end of the simulation with the biased process. 
It is on this equilibrated microscopic sample that the final microscopic accept/reject step will be performed.

Besides its general applicability, an additional advantage of the indirect reconstruction scheme is that we can also use this scheme to efficiently pre-compute quantities that the mM-MCMC requires at the reaction coordinate level: the coefficients in the effective dynamics, and an approximation to the invariant distribution of the reaction coordinate values. Both these quantities can be written as an integral over the sub-manifold of constant reaction coordinate value. The indirect reconstruction then effectively places microscopic samples near this sub-manifold so that these samples can be used for a Monte Carlo approximation of these integrals.

The idea of using an effective dynamics to generate coarse-grained proposals was already proposed in the Coupled Coarse Graining MCMC method, introduced in~\cite{kalligiannaki2012coupled,kalligiannaki2012multilevel}, where large lattice systems with an Ising-type potential energy were sampled. In this specific setting, there are a few natural expressions available that describe the reconstruction distribution, and obtaining a reconstructed `sample' can also be achieved efficiently due to the natural hierarchical nature of a lattice system.
Similarly, a two-level MCMC algorithm is also used in~\cite{efendiev2006preconditioning} as a `pre-conditioner' to increase the microscopic acceptance rate for fluid flows. First, a low-dimensional macroscopic approximation to the high-dimensional fluid simulation is performed, and on acceptance of this macroscopic simulation, the fine-scale fluid simulation is performed afterwards. During the reconstruction step, one can easily reuse some macroscopic basis functions at the microscopic level, making the reconstruction step efficient.

The remainder of this manuscript is organised as follows. In Section~\ref{sec:mM-MCMC}, we briefly summarize the micro-macro Markov chain Monte Carlo method with direct reconstruction that was introduced in \cite{vandecasteele2020direct}. In Section~\ref{sec:maintainingdb}, we introduce the indirect reconstruction scheme and show that a simple replacement of direct reconstruction by indirect reconstruction results in a method that is not reversible and hence is unable to sample from~\eqref{eq:mu}. In Section~\ref{sec:mMMCMCrcs}, we propose a reversible formulation of mM-MCMC with indirect reconstruction on an extended state space and explain the different steps in detail. This section also contains a proof of convergence and a discussion of the influence of the parameters in indirect reconstruction on the efficiency of the mM-MCMC method. In Section~\ref{sec:practical}, we discuss how the indirect reconstruction scheme can be used to pre-compute an approximation to invariant distribution of the reaction coordinate values and the coefficients of the effective dynamics. In Section~\ref{sec:results}, we apply the mM-MCMC scheme with indirect reconstruction to two molecular dynamics cases: an academic three-atom molecule and alanine-dipeptide. In each of these cases, there is a time-scale present between parts of the molecule and we show numerically that mM-MCMC is able to bridge a large part of the time-scale separation. We conclude this manuscript with a summarising discussion and some pointers to future research in Section~\ref{sec:conclusion}.

\section{Micro-macro MCMC with direct reconstruction}\label{sec:mM-MCMC}

In Section \ref{sec:mM-MCMC-direct}, we briefly recap the micro-macro Markov chain Monte Carlo method with direct reconstruction that was introduced in \cite{vandecasteele2020direct}. This method relies on the use of reaction coordinates for the macroscopic description of the system under study.  We therefore first give a short summary and some necessary notation on reaction coordinates in Section~\ref{sec:rc}.

\subsection{A short recap on reaction coordinates}\label{sec:rc}

For the formulation of the mM-MCMC method, we assume  a macroscopic (slow) \emph{reaction coordinate} to be given. A reaction coordinate is a differentiable function from the high-dimensional configuration space $\mathbb{R}^d$ to a lower dimensional space $\mathbb{R}^n$ with $n \ll d$~\cite{stoltz2010free, legoll2010effective}. We denote this function as
\begin{equation}
\xi: \mathbb{R}^d \to \mathbb{R}^n, \ x \mapsto \xi(x) = z,
\end{equation}
where we will also denote by $H \subset \mathbb{R}^n$ the image of $\xi$. Most macroscopic variables in molecular dynamics (e.g., torsion angles, centers of mass between of sets of atoms, \ldots) can be written as a reaction coordinate that is a function of all the positions of the different atoms within the system.

Based on the high-dimensional invariant distribution $\mu$ for the microscopic samples $x$, we can deduce an invariant probability distribution for the reaction coordinate values only. Setting $z=\xi(x)$, the invariant distribution $\mu_0$ of the reaction coordinates reads
\begin{equation} \label{eq:freeenergy}
\mu_0(z) = Z_A^{-1} \exp(-\beta A(z)), \ \ A(z) = -\beta^{-1} \ln\left(\int_{\Sigma(z)} Z_V^{-1} \exp\left(-\beta V(x) \right) \ \norm{\nabla \xi(x)}^{-1} \right) d\sigma_{\Sigma(z)}(x),
\end{equation}
where the integral is taken over the set $\Sigma(z) = \{ x \in \mathbb{R}^d \ | \ \xi(x)=z \}$ of all microscopic samples, conditioned on a given reaction coordinate value $z$. The potential energy function $A(z)$ is called the \emph{(Helmholtz) free energy} of the reaction coordinate.

An important equality that relates integrals over the high-dimensional state space $\mathbb{R}^d$ to integrals over the level sets $\Sigma(z)$ of the reaction coordinate is the \emph{co-area formula}~\cite{stoltz2010free}. Given a function $f \in L^1(\mathbb{R}^d)$, we have the following identity
\begin{equation} \label{eq:coarea}
\int_{\mathbb{R}^d} f(x)  dx = \int_H \int_{\Sigma(z)} f(x) \norm{\nabla \xi(x)}^{-1}  d\sigma_{\Sigma(z)}(x) dz = \int_H \int_{\Sigma(z)} f(x) \delta_{\xi(x)-z} dz .
\end{equation}
in which we formally write, in the sense of measures, $\delta_{\xi(x)-z} = \norm{\nabla \xi(x)}^{-1} d\sigma_{\Sigma(z)}(x)$. 

Besides the invariant distribution of the reaction coordinate values, we can also write down the time-invariant distribution of the microscopic samples, conditioned upon a fixed value $z$ of the reaction coordinate. This distribution is defined on the set $\Sigma(z)$ and is given by
\begin{equation} \label{eq:nu}
\nu(x|z) = \frac{\mu(x)}{\mu_0(z)} \delta_{\xi(x)-z} = \frac{\exp\left(-\beta V(x) \right) \ \norm{\nabla \xi(x)}^{-1}}{\int_{\Sigma(z)} \exp\left(-\beta V(x) \right) \ \norm{\nabla \xi(x)}^{-1}dx }, \ \ x \in \Sigma(z).
\end{equation}
In the last expression, we used the definition of the free energy~\eqref{eq:freeenergy} to rewrite the denominator and the co-area formula~\eqref{eq:coarea} to rewrite the numerator of $\nu$. We will call the distribution $\nu(x|z)$ the \emph{time-invariant direct reconstruction distribution} in this manuscript.

Finally, based on the underlying overdamped Langevin dynamics of the molecular system~\eqref{eq:overdamped_langevin}, one can write down an approximate stochastic evolution equation for the reaction coordinate values by an \emph{effective dynamics}~\cite{legoll2010effective} of the form
\begin{equation}\label{eq:effdyn_intro}
dZ_t = b(Z_t) d + \sqrt{2\beta^{-1}} \sigma(Z_t) dW_t.
\end{equation}
This stochastic process has $\mu_0$, see equation~\eqref{eq:freeenergy}, as invariant distribution and the coefficients $b(z)$ and $\sigma(z)$ are defined by
\begin{align}
b(z) &= \mathbb{E}_\mu[-\nabla V \cdot \nabla \xi + \beta ^{-1} \triangle \xi \ | \ \xi(x) = z] \nonumber  \\   
\sigma(z)^2 &= \mathbb{E}_\mu[\norm{\nabla \xi}^2 \ | \ \xi(x)=z], \label{eq:effdyncoef}
\end{align}
where the expected values are taken over the level set $\Sigma(z)$. In Section~\ref{sec:practical}, we present a straightforward numerical scheme to pre-compute an approximation of these effective dynamics coefficients $b(z)$ and $\sigma(z)$, as well as the free energy $A(z)$, using the indirect reconstruction scheme that is introduced in Section~\ref{sec:mMMCMCrcs}.

\subsection{mM-MCMC with direct reconstruction}\label{sec:mM-MCMC-direct}

To state the complete mM-MCMC method with direct reconstruction, we assume the availability of two ingredients. First, we assume that we can sample an approximation $\bar{\mu}_0(z)$ to the exact invariant probability measure $\mu_0$ of the reaction coordinates, using an MCMC method with a macroscopic transition distribution $q_0(\cdot | \cdot)$. This macroscopic sampling is discussed in Section~\ref{subsubsec:coarseprop}.  Second, we require a reconstruction distribution $\bar{\nu}(x|z)$ of microscopic samples conditioned upon a given reaction coordinate value. The reconstruction step is discussed in Section~\ref{subsubsec:reconstruction}. Both the macroscopic sampling and the reconstruction step involve an accept/reject procedure.  In principle, the choice of $\bar{\mu}_0$ and $\bar{\nu}$ is arbitrary for the mM-MCMC method to converge. However, these choices influence the efficiency of the resulting method. In Section~\ref{subsubsec:transitionkernel}, we give an expression for the  transition kernel of mM-MCMC with direct reconstruction, which we will use in Section~\ref{sec:maintainingdb}.

\subsubsection{Generating a macroscopic proposal}\label{subsubsec:coarseprop}
Let us start from the (given) current microscopic sample $x_n$ on the Markov chain. We generate a new microscopic sample $x_{n+1}$ in two steps. First, we restrict the current microscopic sample to the corresponding reaction coordinate value, i.e., we compute $z_n = \xi(x_n)$.
Next, we sample a new reaction coordinate value $z'$ from the macroscopic transition distribution $q_0(z'|z_n)$. This macroscopic proposal can, for instance, be based on a Brownian motion in the reaction coordinate space or on a gradient descent step using an approximation to the free energy of the reaction coordinate~\eqref{eq:freeenergy}.

To ensure that $z'$ samples the prescribed macroscopic distribution $\bar{\mu}_0$ of the reaction coordinate values, we accept $z'$ with probability
\begin{equation} \label{eq:coarseacceptancerate}
\alpha_{CG}(z'|z_n) = \min\left\{1, \frac{\bar{\mu}_0(z') \ q_0(z_n|z')}{\bar{\mu}_0(z_n) \ q_0(z' |z_n)} \right\}.
\end{equation}
which is the standard Metropolis-Hastings form for the acceptance rate. If $z'$ is rejected, we set $x_{n+1}=x_n$ and generate a new macroscopic proposal; otherwise, we proceed to the reconstruction step.

 \subsubsection{Reconstructing a microscopic sample} \label{subsubsec:reconstruction}
If the reaction coordinate value $z'$ has been accepted, we reconstruct a microscopic sample $x'$ by drawing a sample from the (given) reconstruction distribution $\bar{\nu}(x'|z')$. Then, we accept $x'$ with probability
\begin{equation} \label{eq:microalpha}
\alpha_F(x'|x) = \min\left\{1, \frac{\mu(x') \ \bar{\mu}_0(z_n) \ \bar{\nu}(x_n|z_n)}{\mu(x_n) \ \bar{\mu}_0(z') \ \bar{\nu}(x'|z')}  \right\}.
\end{equation}
This form of the microscopic acceptance has been derived in~\cite{vandecasteele2020direct}. One can easily see that 
\begin{equation} \label{eq:q}
q(x'|x) = \bar{\nu}(x'|\xi(x')) \ \alpha_CG{\xi(x')|\xi(x)} \ q_0(\xi(x')|\xi(x))
\end{equation}
is indeed the transition distribution of generating the microscopic sample $x'$ from $x$. Then, the microscopic acceptance rate~\eqref{eq:microalpha} yields, in fact, the standard Metropolis-Hastings acceptance rate associated with a move from $x$ to $x'$.

On acceptance, define $x_{n+1}=x'$; otherwise set $x_{n+1}=x_n$. With the forms~\eqref{eq:coarseacceptancerate} and~\eqref{eq:microalpha} for the macroscopic and microscopic acceptance rates, the mM-MCMC algorithm with direct reconstruction is guaranteed to sample from $\mu$ correctly~\cite{vandecasteele2020direct}.

\subsubsection{Transition kernel of mM-MCMC with direct reconstruction} \label{subsubsec:transitionkernel}

To finalise the discussion of the mM-MCMC method with direct reconstruction, we introduce the corresponding transition kernel. Based on the definition of the mM-MCMC transition distribution~\eqref{eq:q}, the total probability of transitioning from $x$ to $x'$ is
\begin{equation} \label{eq:mMtransitionkernel}
\mathcal{K}_{mM}(x'|x) = \begin{cases}
& \alpha_F(x'|x) \ \bar{\nu}(x'|\xi(x')) \ \alpha_{CG}(\xi(x')|\xi(x)) \ q_0(\xi(x')|\xi(x)), \ \ x' \neq x \\
\int_{\mathbb{R}^d} &\alpha_F(y|x) \ \bar{\nu}(y|\xi(y)) \ \alpha_{CG}(\xi(y)|\xi(x)) \ q_0(\xi(y)|\xi(x)) \ dy, \ \ x' = x 
\end{cases}.
\end{equation}

\section{An indirect reconstruction scheme using biased simulation} \label{sec:maintainingdb}

An important ingredient for the micro-macro Markov chain Monte Carlo algorithm with direct reconstruction is the choice of reconstruction distribution $\bar{\nu}$. This reconstruction distribution should be easy to sample and should result in a high microscopic acceptance rate. 
However, finding a reconstruction distribution that satisfies these criteria may be highly non-trivial and application-dependent. The time-invariant reconstruction distribution $\nu$ is in principle always available, but is usually hard to sample from. Indeed, sampling from the reconstruction distribution
\[
\nu(x'|  z') \propto \exp\left(-\beta V(x')\right) \ \norm{\nabla \xi(x')}^{-1}, \ \ x \in \Sigma(z')
\]
is often harder than sampling the target distribution $\mu$, because the domain $\Sigma(z')$ can have an irregular form, depending on the form of $\xi$. In this section, we therefore propose a more general \emph{indirect reconstruction scheme} that reconstructs a microscopic sample close to the sub-manifold $\Sigma(z')$, but only approximately corresponds to the desired reaction coordinate value.

In Section~\ref{subsec:matching}, we discuss the biased stochastic process on which the indirect reconstruction is based. The microscopic sample $x'$, obtained by the indirect reconstruction scheme, does not generally lie on the given sub-manifold of constant reaction coordinate value $\Sigma(z')$. We show in Section~\ref{subsec:ensuringdb} that, as a consequence, a simple extension of mM-MCMC with direct reconstruction to the indirect reconstruction setting does not result in a reversible scheme. We solve the issue of reversibility in Section~\ref{subsubsec:matchingkernel} by making extending the state space of the mM-MCMC method.

\subsection{A straightforward scheme for indirect reconstruction} \label{subsec:matching}
Consider a reaction coordinate value $z'$ obtained by the macroscopic proposal step at the macroscopic level. To reconstruct a microscopic sample associated with $z'$ via indirect reconstruction, we perform a time integration of the biased stochastic process
\begin{equation} \label{eq:constrainedsimulation}
dX  = -\nabla \left(V(X) + \frac{\lambda}{2}\norm{\xi(X)-z'}^2\right) dt + \sqrt{2\beta^{-1}} dW,
\end{equation}
that pulls the reaction coordinate value of $X$ towards $z'$. This biased process has an invariant distribution of the form 
\begin{equation} \label{eq:constrainedinvariant}
\nu_\lambda(x'; z') = \left(\frac{\lambda \beta}{2\pi}\right)^{n/2}  \frac{1}{N_\lambda(z')}\exp \left(-\beta V(x')\right) \ \exp\left( -\frac{\lambda \beta}{2}\norm{\xi(x')-z'}^2 \right),
\end{equation}
which we will also call the \emph{indirect reconstruction distribution} in this manuscript. The value $N_\lambda(z')$ is a normalization constant depending on $z'$ and $\lambda$. Note that we introduced the notation  $\nu_\lambda(x'; z')$ for the indirect reconstruction distribution because the microscopic sample $x' \in \mathbb{R}^d$ does not necessarily have reaction coordinate $z'$; instead $z'$ merely acts as a parameter in the invariant distribution of the process~\eqref{eq:constrainedinvariant}.
The factor $\exp \left(-\beta V(x')\right)$ in \eqref{eq:constrainedinvariant} corresponds to the invariant distribution~\eqref{eq:mu} of the microscopic process~\eqref{eq:overdamped_langevin}, while the extra term $\frac{\lambda}{2}\norm{\xi(x')-z'}^2$ in the potential energy ensures that the measure $\nu_\lambda(x'; z')$ will concentrate on values of $x'$ for which $\xi(x')$ is close to $z'$, provided $\lambda$ is large enough.

Before giving a numerical discretization of the biased stochastic process, we first give an expression of the normalization constant $N_\lambda(z')$. Integrating~\eqref{eq:constrainedinvariant} using the co-area formula gives
\begin{align}
&\int_{\mathbb{R}^d} \left(\frac{\lambda \beta}{2\pi}\right)^{n/2} \frac{1}{N_\lambda(z')} \exp(-\beta V(x')) \exp\left(-\frac{\beta \lambda}{2} \norm{\xi(x')-z'}^2\right) \ dx'  \nonumber \\
&= \frac{Z_V}{Z_A} \int_{H} \left(\frac{\lambda \beta}{2\pi}\right)^{n/2}  \frac{\exp\left(-\frac{\lambda \beta}{2}(u-z')^2\right) \exp\left(-\beta A(u))\right)}{N_\lambda(z')}  \int_{\Sigma(u)} \frac{Z_A}{Z_V}\frac{\exp\left(-\beta (V(x') - A(u))\right)}{\norm{\nabla \xi(x')}} d\sigma_{\Sigma(u)}(x') \ du. \nonumber
\end{align}
The second integral over $\Sigma(u)$ is precisely the time-invariant direct reconstruction distribution $\nu(x' | u)$, see equation~\eqref{eq:nu}, which integrates to $1$. We thus obtain
\begin{equation*}
\int_{\mathbb{R}^d} \nu_\lambda(x'; z') dx' = \frac{Z_V}{Z_A} \int_{H} \left(\frac{\lambda \beta}{2\pi}\right)^{n/2}  \frac{\exp\left(-\frac{\lambda \beta}{2}(u-z')^2\right) \exp\left(-\beta A(u))\right)}{N_\lambda(z')} du = 1. \nonumber
\end{equation*}
Rewriting the expression above, the formula for the normalization constant $N_\lambda(z')$ reads
\begin{equation} \label{eq:normalizationconstant}
N_\lambda(z') = \frac{Z_V}{Z_A} \int_{H} \left(\frac{\lambda \beta}{2\pi}\right)^{n/2}  \exp\left(-\frac{\lambda \beta}{2}(u-z')^2\right)  \exp\left(-\beta A(u))\right) du.
\end{equation}
The expression for $N_\lambda(z')$ is a convolution of a Gaussian form with the invariant distribution of the reaction coordinates $\mu_0(z')$. This convolution operation is also called a scaled Weierstrass transform~\cite{hirschman2012convolution} and can be viewed as a filter on $\mu_0$. In Section~\ref{sec:mMMCMCrcs}, it will become clear that the value of $N_\lambda(z')$ needs to be known in the mM-MCMC method with indirection reconstruction. In Section~\ref{sec:practical}, we give a simple numerical scheme to compute the normalization constant $N_\lambda$ based on the indirect reconstruction scheme itself.

Starting with the previous microscopic sample obtained by the mM-MCMC method $x_n = x_{n,0}$, we generate a sequence of microscopic samples $x_{n,k}, k = 0, \dots, K$ by time-stepping the biased dynamics~\eqref{eq:constrainedsimulation}, using the Euler-Maruyama discretization
\begin{equation} \label{eq:constrainedsimulation_discrete}
x_{n, k+1} = x_{n,k} - \nabla V(x_{n,k})\delta t - \lambda (\xi(x_{n,k})-z') \nabla \xi(x_{n,k}) \delta t + \sqrt{2\beta^{-1} \delta t} \ \eta_{n, k} , \ \ \eta_{n,k} \sim \mathcal{N}(0, 1),
\end{equation}
followed by an accept/reject step to ensure the microscopic samples $x_{n,k}$ indeed sample the indirect reconstruction distribution $\nu_\lambda(\cdot; \ z')$. That is, we accept $x_{n,k+1}$ from $x_{n,k}$ with probability
\[
\min\left\{1, \frac{\nu_\lambda(x_{n,k+1}; z') \ q_{EM}(x_{n,k} | x_{n,k+1})}{\nu_\lambda(x_{n,k}; z') \ q_{EM}(x_{n,k+1} | x_{n,k})}  \right\},
\]
where $q_{EM}$ is the Euler-Maruyama transition distribution based on the discretization~\eqref{eq:constrainedsimulation_discrete}.
This sampling scheme is an example of the MALA method (Metropolis adjusted Langevin algorithm)~\cite{roberts1996exponential}. We perform $K$ time steps of the biased simulation process to overcome the burn-in period associated with the inconsistency of the initial condition $x_n$ with the desired reaction coordinate value $z'$. After these $K$ steps, we propose the microscopic sample $x_{n, K}$ to be the reconstructed microscopic sample, i.e.,  $x' = x_{n, K}$. This final step of the biased simulation is then followed by the microscopic accept/reject step to decide on the acceptance of $x'$. We give more details on the acceptance criterion and the complete mM-MCMC with indirect reconstruction in Section~\ref{sec:mMMCMCrcs}. 

The efficiency and accuracy of the indirect reconstruction scheme depend on the choice of parameters $\lambda$ and $K$ and the time step $\delta t$ used to simulate~\eqref{eq:constrainedsimulation_discrete}. In Section~\ref{subsec:optimalparameters}, we provide heuristics to determine a suitable set of parameters  for the indirect reconstruction scheme.

\subsection{A simple formulation of mM-MCMC with indirect reconstruction is not reversible} \label{subsec:ensuringdb}

With the indirect reconstruction scheme, a reconstructed microscopic sample $x'$ does generally not correspond to the proposed reaction coordinate value $z'$ at the macroscopic level, i.e., we no longer have the relation $\xi(x')=z'$ after reconstructing $x'$ from $z'$. As a consequence, there is a fundamental asymmetry between restriction and reconstruction. On one hand, each microscopic sample corresponds to only one reaction coordinate value. On the other hand, however, after indirect reconstruction, the reaction coordinate value corresponding to the reconstructed microscopic state in $\mathbb{R}^d$ can, in principle, be any reaction coordinate value. Therefore, if we simply plug in the indirect reconstruction distribution in the mM-MCMC method with direct reconstruction, i.e., if we use the reconstruction distribution $\bar{\nu} = \nu_\lambda$ in the microscopic acceptance rate~\eqref{eq:microalpha}, one can expect that the resulting method stops being reversible. In this Section, we show this intuitive result mathematically. To restore reversibility, we propose a modification to the mM-MCMC method with indirect reconstruction in Section~\ref{subsubsec:matchingkernel}.

Indeed, suppose that we replace the reconstruction distribution $\bar{\nu}$ by the indirect reconstruction distribution $\nu_\lambda$ in the transition kernel $\mathcal{K}_{mM}$~\eqref{eq:mMtransitionkernel} of mM-MCMC with direct reconstruction. The resulting expression for the transition kernel is then not only a function of the microscopic samples $x$ and $x'$, but also depends on the reaction coordinate value $z'$ sampled at the macroscopic level. Thus, to obtain a closed expression for the mM-MCMC transition kernel with indirect reconstruction only in $x$ and $x'$, we integrate this transition kernel $\mathcal{K}_{mM}$~\eqref{eq:mMtransitionkernel} over all possible values of the reaction coordinate $z'$. This integrated transition kernel with the indirect reconstruction method reads
\begin{equation}\label{eq:kernel_direct}
\mathcal{K}(x'|x) =\begin{cases*} \begin{aligned} &\int_H \alpha_F(x'| x, z') \ \nu_\lambda(x'; z') \  \alpha_{CG}(z'|\xi(x)) \ q_0(z'|\xi(x)) dz' \ \ x' \neq x \\
& 1- \int_{\mathbb{R}^d} \int_H \alpha_F(y| x, z') \ \nu_\lambda(y; z') \ \alpha_{CG}(z'|\xi(x)) \ q_0(z'|\xi(x)) dz' \ dy, \ \ x' = x,
\end{aligned}
\end{cases*}
\end{equation}
where the microscopic acceptance rate $\alpha_F(x'|x,z')$ has the same form as for direct reconstruction, but we need the extra dependency on $z'$ since $\xi(x') \neq z'$, i.e., 
\[
\alpha_F(x'|x,z') = \left\{1, \frac{\mu(x') \ \bar{\mu}_0(\xi(x)) \ \nu_\lambda(x; \xi(x))}{\mu(x) \ \bar{\mu}_0(z') \ \nu_\lambda(x'; z')} \right\}.
\]

Algorithmically speaking, the mM-MCMC algorithm with indirect reconstruction given by the above transition kernel $\mathcal{K}$ is identical to the direct reconstruction method that we summarised in Section~\ref{sec:mM-MCMC-direct}. Only the reconstruction step is different. The transition kernel~\eqref{eq:kernel_direct}, however, fails so sample the microscopic Gibbs measure $\mu(dx)$ because it does not satisfy detailed balance. Indeed, if we write out the detailed balance condition for $x' \neq x$, we find
\begin{equation*}
\begin{aligned}
\mathcal{K}(x'|x) \ \mu(x) &= \int_H \min\left\{1,\frac{\mu(x') \ \bar{\mu}_0(\xi(x)) \ \nu_\lambda(x; \xi(x))}{\mu(x) \ \bar{\mu}_0(z') \ \nu_\lambda(x'; z')} \right\} \ \nu_\lambda (x'; z') \\ &\times \min\left\{1,\frac{\bar{\mu}_0(z') \ q_0(\xi(x) | z')}{\bar{\mu}_0(\xi(x)) \ q_0(z'|\xi(x))}  \right\} \ q_0(z'|\xi(x)) \mu(x) dz' \\
&= \int_H \min\left\{\mu(x) \ \bar{\mu}_0(z') \ \nu_\lambda(x'; z'), \mu(x') \ \bar{\mu}_0(\xi(x)) \ \nu_\lambda(x; \xi(x)) \right\} \\
&\times \min\left\{\frac{q_0(z'|\xi(x))}{\bar{\mu}_0(z')}, \frac{q_0(\xi(x)|z')}{\bar{\mu}_0(\xi(x))} \right\} dz'.
\end{aligned}
\end{equation*}
The above expression is not symmetric in $x$ and $x'$ because $\nu_\lambda(x'; z')$ is integrated over all possible $z$, while $\nu_\lambda(x; \xi(x))$ is independent of $z'$. Therefore, writing the mM-MCMC scheme with indirect reconstruction only in terms of the microscopic samples does not result in a reversible Markov chain. This above observation motivates the need to formulate the mM-MCMC method on an extended state space that includes both the microscopic samples and the reaction coordinate values as decoupled variables, as will be done in the next Section. 

\section{Micro-macro MCMC method indirect reconstruction} \label{sec:mMMCMCrcs}
Having introduced the indirect reconstruction scheme and having shown its renders the mM-MCMC algorithm of Section~\ref{sec:mM-MCMC-direct} irreversible, we now state a modified mM-MCMC algorithm with indirect reconstruction on an extended state space, so that the resulting method becomes reversible. 
 The extended state space consists of tuples $(x, z)$ of a microscopic sample $x$ and a reaction coordinate value $z$, and is given by $\mathbb{R}^d \times H$, with $H$ the image of the reaction coordinate $\xi$.

In Section~\ref{subsubsec:matchingkernel}, we explain the different steps in the mM-MCMC method with indirect reconstruction and we also summarize the complete mM-MCMC algorithm with indirect reconstruction in Algorithm~\ref{algo:mM-MCMCind}.
In Section~\ref{subsubsec:convergence4matching}, we give an expression of the time-invariant distribution of the mM-MCMC method on the extended state space, where we explicitly show that the marginal time-invariant distribution of the microscopic samples is the Gibbs distribution~\eqref{eq:mu} that we set out to sample from. We also state and prove that the mM-MCMC scheme with indirect reconstruction indeed converges to the correct probability distribution and that the method is ergodic. Finally, the indirect reconstruction method contains a number of additional numerical parameters that have an impact on the overall efficiency of mM-MCMC. Section~\ref{subsec:optimalparameters} contains a discussion on an a priori suitable choice for these parameter values.

\subsection{mM-MCMC with indirect reconstruction on an extended state space} \label{subsubsec:matchingkernel}

Given the current tuple $(x_n, z_n) \in \mathbb{R}^d \times H$ of the Markov chain on the extended state space, we again propose a new tuple $(x_{n+1}, z_{n+1})$ on the Markov chain by (i) generating a new macroscopic proposal (Section~\ref{sec:mM-MCMC-ind-prop}); and (ii) reconstructing a microscopic sample (Section~\ref{sec:mM-MCMC-ind-rec}).
Before detailing these steps, we point out that the restriction step is eliminated here because of the extension of the state space. That is, we do not need to compute $\xi(x_n)$ since we propose new macroscopic values based on the associated reaction coordinate value $z_n$, which is already available as an independent variable the tuple $(x_n, z_n)$. 

\subsubsection{Generating a macroscopic proposal}\label{sec:mM-MCMC-ind-prop}

In the first step, we propose a new reaction coordinate value $z'$ based on $z_n$, using a macroscopic transition distribution $q_0(z'|z_{n})$. To ensure that $z'$ is sampled according to the approximate macroscopic distribution $\bar{\mu}_0$, we accept $z'$ with probability $\alpha_{CG}$~\eqref{eq:coarseacceptancerate}.

If the macroscopic proposal $z'$ is accepted, we proceed to the reconstruction step, otherwise, we immediately define $(x_{n+1}, z_{n+1}) = (x_n, z_n)$ and repeat this step. Note that this step is identical to the macroscopic proposal step in the direct reconstruction algorithm in Section~\ref{subsubsec:coarseprop}, since the reconstruction has no impact on the macroscopic proposal moves.

\subsubsection{Reconstructing a microscopic sample}\label{sec:mM-MCMC-ind-rec}
In the second step, we reconstruct a new microscopic sample $x'$ from the proposed reaction coordinate value $z'$ using the indirect reconstruction step discussed in Section~\ref{sec:mM-MCMC-ind-prop}. That is, we simulate $K$ MALA steps~\eqref{eq:constrainedsimulation_discrete} with time step $\delta t$ and with initial condition $x_n$. Afterwards, we decide whether to accept $x'$ using an acceptance criterion.

If we use the same form of the microscopic acceptance probability as defined in the direct reconstruction algorithm~\eqref{eq:microalpha} but plug in the indirect reconstruction distribution $\nu_\lambda$, the microscopic acceptance probability would read
\[
\alpha_F(x'|x, z', z) = \min\left\{ 1, \frac{\mu(x') \ \bar{\mu}_0(z) \ \exp\left(-\beta V(x)\right) \ \exp\left(-\frac{\beta \lambda}{2}\left((\xi(x)-z)^2\right)\right) \ N_\lambda(z')}{\mu(x) \ \bar{\mu}_0(z') \ \exp\left(-\beta V(x')\right) \ \exp\left(-\frac{\beta \lambda}{2}\left((\xi(x')-z')^2\right) \right) \ N_\lambda(z')}  \right\}.
\]
One can see that when $\lambda$ is large, the ratio of the two sharp Gaussian factors could significantly lower the average microscopic acceptance rate when $\xi(x)$ lies farther from $z$ than $\xi(x')$ lies from $z'$. In this case, the complete mM-MCMC method with indirect reconstruction would become inefficient since we would then reject many reconstructed microscopic samples that are computationally expensive to obtain.

To overcome this computational issue, we use a modified microscopic acceptance probability where we remove the ratio of the sharp Gaussians. That is, we accept the reconstructed microscopic sample $x'$ with the modified probability
\begin{align}
\alpha_F(x'| x, z', z) &= \min\left\{1, \frac{\mu(x') \ \bar{\mu}_0(z) \ \exp\left(-\beta V(x)\right) \ N_\lambda(z')}{\mu(x) \ \bar{\mu}_0(z') \ \exp\left(-\beta V(x')\right) \ N_\lambda(z')} \right\} \nonumber \\ &= \min \left\{1, \frac{\exp(-\beta V(x')) \ \bar{\mu}_0(z) \ \exp\left(-\beta V(x)\right)   N_\lambda(z')}{\exp(-\beta V(x)) \ \bar{\mu}_0(z') \ \exp\left(-\beta V(x')\right) \ N_\lambda(z)} \right\} \nonumber \\
&= \min \left\{ 1, \ \frac{\bar{\mu}_0(z)}{\bar{\mu}_0(z')} \frac{N_\lambda(z')}{N_\lambda(z)} \right\}. \label{eq:matchingacceptance}
\end{align}
We will show in Section~\ref{subsubsec:convergence4matching} that the mM-MCMC method with this acceptance probability is indeed reversible with respect to $\mu$ and that only the marginal invariant distribution of the reaction coordinate values is changed. We will also show in the numerical experiments that the average microscopic acceptance rate is close the $1$ as desired.

Finally, on acceptance we define $(x_{n+1}, z_{n+1}) = (x', z')$, otherwise we set $(x_{n+1}, z_{n+1}) = (x_n, z_n)$ and return to the macroscopic proposal step.

\begin{remark}
In case the macroscopic invariant distribution is the marginal time-invariant distribution of the reaction coordinate of the full dynamics, i.e., $\bar{\mu}_0 = \mu_0$, the expression for the microscopic acceptance rate converges to $1$ as $\lambda$ increases to $\infty$. Indeed, we have the limit 
\[
\lim\limits_{\lambda \to \infty} N_\lambda(z') = \frac{Z_V}{Z_A}\exp(-\beta A(z')) = Z_V \mu_0(z').
\]
Hence, if we drive $\lambda$ to infinite, the microscopic acceptance rate~\eqref{eq:matchingacceptance} becomes
\[
\begin{aligned}
\lim\limits_{\lambda \to \infty} \min\left\{ 1, \ \frac{\mu_0(z)}{\mu_0(z')} \frac{N_\lambda(z')}{N_\lambda(z)}\right\} &= \min\left\{1, \lim\limits_{\lambda \to \infty}\frac{\mu_0(z)}{\mu_0(z')} \frac{N_\lambda(z')}{N_\lambda(z)} \right\} \\
&= \min\left\{1, \frac{\mu_0(z) \ \mu_0(z')}{\mu_0(z') \ \mu_0(z)} \right\} = 1,
\end{aligned}
\]
in point-wise sense. We can interchange the limit and the minimum in the first equality due to the fact that both functions inside the minimum are continuous in $z$ and $z'$, and so is the minimum of both functions.
\end{remark}

\subsubsection{Complete algorithm\label{sec:alg2}}
The complete mM-MCMC algorithm with indirect reconstruction is depicted in Algorithm~\ref{algo:mM-MCMCind}. 
\begin{algorithm}[h!] 
	\begin{flushleft}
		Given parameters $\lambda > 0$, $\delta t>0$, $K \in \mathbb{N}$, microscopic samples $x_n, \ n = 1, 2, \dots$ and reaction coordinate values $z_n, \ n=1, 2, \dots$, with $z_0 = \xi(x_0)$.
		
		\vspace{4mm}
		(i) \textbf{Macroscopic Poposal}: 
		
		\begin{itemize}
			\item Generate a reaction coordinate proposal $z' \sim q_0 ( \cdot | z_n)$.
			
			\item Accept the macroscopic proposal with probability
			\[
			\alpha_{CG}(z' |  z_n) = \min \left\{1, \frac{\bar{\mu}_0(z') \ q_0(z_n | z')}{\bar{\mu}_0(z_n) \ q_0(z' | z_n)} \right\}
			\]
			
			\item On acceptance, proceed to step (ii), otherwise set $(x_{n+1}, z_{n+1}) = (x_n, z_n)$ and repeat step (i).
		\end{itemize}
		
		\vspace{4mm}
		(ii) \textbf{Reconstruction}: 
		\begin{itemize}
		\item Perform $K$ MALA steps of the biased simulation~\eqref{eq:constrainedsimulation_discrete} with step size $\delta t$ and initial condition $x_{n,0} = x_n$
		\[
		x_{n,k+1} = x_{n,k} - \nabla \left(V(x_{n,k}) + \frac{\lambda}{2} \norm{\xi(x_{n,k})-z'}^2\right)\delta t + \sqrt{2\beta^{-1} \delta t} \ \eta_{n,k} , \ \ \ k = 0, \dots,K-1, \ \eta_{n,k} \sim \mathcal{N}(0,1)
		\]
		and define the microscopic sample $x' = x_{n,K}$.  
		
		\item Accept the microscopic sample with probability
		\[
		\alpha_F(x'| x_n, z', z_n) = \min \left\{ 1, \ \frac{\bar{\mu}_0(z_n)}{\bar{\mu}_0(z')} \frac{N_\lambda(z')}{N_\lambda(z_n)}  \right\}.
		\]
		
		\item Upon acceptance, set $(x_{n+1}, z_{n+1}) = (x', z')$ and return to step (i) for the next microscopic sample, otherwise set $(x_{n+1}, z_{n+1}) = (x_n, z_n)$ and generate a new reaction coordinate value in step (i).
		\end{itemize}
	\end{flushleft}
	\caption{The mM-MCMC method with indirect reconstruction.}
	\label{algo:mM-MCMCind}
\end{algorithm}

 \subsection{Convergence result} \label{subsubsec:convergence4matching}
 In this Section, we show that the mM-MCMC scheme with indirect reconstruction has a unique time-invariant distribution and that the method is ergodic. First, we given an expression for the transition kernel of mM-MCMC with indirect reconstruction. This transition kernel, defined on the extended state space $\mathbb{R}^d \times H$, reads
 \begin{equation} \label{eq:matchingtransitionkernel}
 \mathcal{K}^*_{mM}(x', z' | x, z) = \begin{cases}
 \begin{aligned}
 &\alpha_F(x'|x, z', z) \ \nu_\lambda(x'; z') \ \alpha_{CG}(z'|z) \ q_0(z'|z) \ \ (x', z') \neq (x, z) \\
 &1- \int_{\mathbb{R}^d} \int_H \alpha_F(y|x, u, z) \ \nu_\lambda(y; u) \ \alpha_{CG}(u|z) \ q_0(u|z) du \ dy \ \  (x', z') = (x, z).
 \end{aligned}
 \end{cases}
 \end{equation}
 
We show in the proof of Theorem~\ref{thm:convergence_matching} that there is a unique time-invariant distribution associated to the transition kernel $\mathcal{K}_{mM}^*$ mM-MCMC with indirect reconstruction. We additionally show that the mM-MCMC method converges to this distribution under some mild conditions on the macroscopic transition distribution $q_0$ and the approximate macroscopic distribution $\bar{\mu}_0$.
 
\begin{theorem} \label{thm:convergence_matching}
	The mM-MCMC transition kernel $\mathcal{K}_{mM}^*$ has a unique stationary probability distribution of the form
	\begin{equation} \label{eq:eta}
		\eta(x, z)  = \sqrt{\frac{\beta \lambda}{2\pi}}  \ \exp\left(-\frac{\beta \lambda}{2} \left(z-\xi(x)\right)^2\right) \ \mu(x).
	\end{equation}
	Furthermore, for every macroscopic transition distribution $q_0$ that is strictly positive on $H\times H$ and that is not identical to the exact, time-discrete, transition distribution of the effective dynamics~\eqref{eq:effdyn_intro}, and for every approximate macroscopic distribution $\bar{\mu}_0$ with $\text{supp}(\bar{\mu}_0) = H$,
	\begin{itemize}
		\item[(i)] the transition probability kernel $\mathcal{K}_{mM}^*$~\eqref{eq:matchingtransitionkernel} satisfies the detailed balance condition with target measure $\eta(x, z)$;
		
		\item[(ii)] the chain $(x_n, z_n)$ is $\eta-$irreducible;
		
		\item[(iii)] the chain $(x_n, z_n)$ is aperiodic.
	\end{itemize}
\end{theorem}
 Note that we can write~\eqref{eq:eta} as 
$\eta(x, z)  = \eta(z|x)\mu(x)$,
 to view $\mu(x)$ as the marginal distribution of the microscopic samples of $\eta(x, z)$. This implies we are indeed sampling the microscopic Gibbs measure~\eqref{eq:mu} that we set out to sample from. The conditional distribution $\eta(z|x)=\sqrt{\frac{\beta \lambda}{2\pi}} \exp\left(-\frac{\beta \lambda}{2} \left(z-\xi(x)\right)^2\right)$ is the time-invariant conditional distribution of reaction coordinate values given a microscopic sample, and is given by a Gaussian distribution with mean $\xi(x)$ and variance $1/(\lambda \beta)$. This identity indicates that the reaction coordinate value of $x$ lies close to the sampled value $z$ at the macroscopic level, which is a logical consequence of the form of the indirect reconstruction distribution. Therefore, in the limit of $\lambda$ increasing to infinity, we recover the mM-MCMC algorithm with direct reconstruction.

Additionally, if we decompose $\eta(x, z) = \eta(z) \eta(x|z)$, i.e., as the marginal distribution of the reaction coordinate values at the macroscopic level multiplied by the conditional distribution of the microscopic samples, one can easily see that 
\begin{equation*}
\begin{aligned}
&\eta(z) = N_\lambda(z), \\
&\eta(x|z) = \nu_\lambda(x; z) \\
\end{aligned}
\end{equation*}
so that the marginal distribution of the reaction coordinate values alone is not the exact marginal distribution $\mu_0$~\eqref{eq:freeenergy}, but the filtered distribution $N_0$ obtained by the Weierstrass transform~\eqref{eq:normalizationconstant}.

\begin{proof}
	We first show that the detailed balance condition holds because this is a sufficient condition for $\eta$ to be a stationary distribution of $\mathcal{K}_{mM}^*$. The combination of statements (ii) and (iii) is then sufficient for $\eta$ to be unique and to show that the mM-MCMC scheme with indirect reconstruction is ergodic and hence converges to $\eta(x,z)$ in total variation norm.
	
	(i) If we denote the total normalization constant $Z_V \sqrt{\frac{2\pi}{\beta\lambda}}$ of $\eta(x, z)$ by $Z_\lambda$, the detailed balance condition reads
	\begin{equation*}
	\begin{aligned}
	\mathcal{K}_{mM}^*(x', z'|x,z) \ \eta(x,z) &= Z_\lambda^{-1}\min\left\{1, \frac{\bar{\mu}_0(z) N_\lambda(z')}{\bar{\mu}_0(z') N_\lambda(z)} \right\} \frac{ \exp\left(-\beta V(x')\right) \exp\left(-\frac{\beta\lambda}{2} (\xi(x')-z')^2\right)}{N_\lambda(z')} \\
	&\times \min\left\{1,\frac{\bar{\mu}_0(z') q_0(z|z')}{\bar{\mu}_0(z) q_0(z'|z)} \right\} q_0(z'|z)  \exp\left(-\beta V(x)\right) \exp\left(-\frac{\beta\lambda}{2} (\xi(x)-z)^2\right) \\
	&= Z_\lambda^{-1} \min\left\{\frac{N_\lambda(z')^{-1}}{ \bar{\mu}_0(z)}, \frac{N_\lambda(z)^{-1} }{\bar{\mu}_0(z')} \right\}\exp\left(-\frac{\lambda \beta}{2} \left( (\xi(x')-z')^2 + (\xi(x)-z)^2\right) \right)  \\
	&\times \exp\left(-\beta V(x) -\beta V(x')\right)  \min\left\{ q_0(z'|z)\bar{\mu}_0(z), q_0(z|z')\bar{\mu}_0(z') \right\} \\
	\end{aligned}
	\end{equation*}
	The above expression is indeed symmetric when interchanging $(x, z)$ and $(x',z')$ and hence detailed balance is satisfied. As a consequence, this fact implies the first claim of the theorem.
	
	\vspace{1mm}
	(ii) To prove $\eta-$irreducibility, take a set measurable set $A \subset \mathbb{R}^d \times H$ with $\eta(A) > 0$ and take an $(x,z) \in \mathbb{R}^d \times H$.
	We then have that 
	\[
	\mathcal{K}^*_{mM}(A|  x, z) \geq \int_{A\backslash \{(x,z)\}}  \alpha_F(x'|x, z',z) \ \nu_\lambda(x' ; z') \ \alpha_{CG}(z'| z) \ q_0(z'| z)\ dx' \ dz'.
	\]
	By assumption, $\nu_\lambda(x';  z'), q_0(z' | z), \alpha_{CG}(z'| z)$ and $\alpha_F(x'| x, z', z)$ are positive everywhere, proving that $\mathcal{K}^*_{mM}(A| x,z) > 0$ for all $(x,z) \in \mathbb{R}^d\times H$.
	
	\vspace{1mm}
	(iii) Finally, we establish aperiodicity of the Markov chain by deriving a contradiction. Assume that for all $(x,z) \in \mathbb{R}^d\times H$ $\mathcal{K}^*_{mM}(\{(x,z)\}| x,z) = 0$. Then, the equality
	\[
	\int_{\mathbb{R}^d} \int_H   \alpha_F(x'| x, z', z) \ \nu_\lambda(x'; z') \ \alpha_{CG}(z'| z) \ q_0(z' | z) dz' \  dx' = 1,
	\]
		 implies the acceptance probabilities $\alpha_F(x'|x,z', z)$ and $\alpha_{CG}(z'|z)$ should be $1$ everywhere because $\nu_\lambda(x';  z')$ and $q_0(z'|z)$ are strictly positive. These two proposal distributions are strictly positive on the whole domain by assumption, implying that the proposal distributions $\nu_\lambda(x';  z')$ and $q_0(z'|z)$ sample from the correct invariant distribution $\eta$ without rejections. This is not the case because the macroscopic transition distribution $q_0$ is not identical to the exact, time-discrete, transition distribution of the effective dynamics~\eqref{eq:effdyn_intro}, resulting in a contradiction.
\end{proof}

\subsection{Optimal choice of parameters for indirect reconstruction} \label{subsec:optimalparameters}
There is still some freedom regarding the choice of parameters in the indirect reconstruction scheme:
\begin{itemize}
	\item What are optimal values for $\lambda$ and $\delta t$ to maximize the efficiency of indirect reconstruction?

	\item How many biased steps $K$ are required to overcome the effect of the burn-in during the biased simulation and to form a sample of~\eqref{eq:constrainedinvariant}?
\end{itemize}

We give more context to these questions in this section, and propose some heuristics for an efficient implementation.

\paragraph{Optimal values of $\lambda$ and $\delta t$}
There are two interesting limit cases for the choice of $\lambda$. In the limit of $\lambda$ decreasing to $0$, the Gaussian factor in the reconstruction distribution~\eqref{eq:constrainedinvariant} has no impact, and the biased simulation is then nothing more than the microscopic MALA method with invariant distribution $\mu$. We then make no use at all of the sampled reaction coordinate value $z'$ at the macroscopic level, and we may then expect to obtain no efficiency gain. 
In this limit, the effect of the bias term  is smaller than that of the stiff modes in the system, annihilating the potential speed-up generated by the reaction coordinates at the macroscopic level.

In the other limit, when $\lambda$ is much larger than the stiffest modes in the potential energy function $V$, we need many biased steps for the reaction coordinate to become close to the sampled value $z'$ due to stiffness of the resulting biased dynamics. Indeed, we must choose the microscopic time step $\delta t$ on the order of $1/\lambda$ due to the stability restrictions on the Euler-Maruyama scheme~\eqref{eq:constrainedsimulation_discrete}, so that one needs many biased time steps to achieve a thorough mixing in the microscopic state space. In this case, the indirect reconstruction scheme is also not efficient anymore.

This effect is also visible in another manner. When $\lambda$ is infinite, we have in the limit on $\mathbb{R}^d$
\begin{align*}
\exp\left(-\beta V(x')\right) \sqrt{\frac{\lambda \beta}{ 2 \pi}} \exp\left(-\frac{\lambda \beta}{2} \norm{\xi(x')-z'}^2\right)dx &\to \exp\left(-\beta V(x')\right) \delta_{\xi(x')-z'}(dx')  \\ &= \exp\left(-\beta V(x')\right) \norm{\nabla \xi(x')}^{-1} d\sigma_{\Sigma(z')}(x'),
\end{align*}
in point-wise sense because the Gaussian factor $\exp\left(-\frac{\lambda \beta}{2} \norm{\xi(x')-z'}^2\right)$ converges to $\delta_{\xi(x')-z'}$. We are then effectively sampling the direct reconstruction distribution $\nu(\cdot| \ z)$ on the set $\Sigma(z)$, which we know to be inefficient in many molecular applications.

Making a trade-off between these two effects, we propose to choose $\lambda$ close to the stiffest mode in the potential energy $V$. This choice should ensure a good mixing in all components of the system, while the microscopic samples $x_{n,k}$ will quickly have a reaction coordinate value close to $z'$. In addition, we choose $\delta t$ close to $1/\lambda$ so that we quickly converge to $z'$, while maintaining a stable biased simulation scheme. We will show numerically in Section~\ref{sec:results} that there is usually a range of values for $\lambda$ that give a maximal, or close to maximal, efficiency gain.

\paragraph{Number of biased steps K}
Deciding on a good number of biased time steps $K$ is a theoretically hard problem, because it depends on the step size $\delta t$ and the strength of the basing potential $\lambda$. Indeed, for any finite number of time steps $K$, the microscopic sample obtained by indirect reconstruction is not in equilibrium with respect to the indirect reconstruction distribution $\nu_\lambda(\cdot; z')$. Hence, there is always a bias present due to indirect reconstruction, and this bias may heavily depend on $\lambda$ and $\delta t$. However, we can make this bias arbitrarily small by making $K$ large enough. Usually, we observe numerically that $5$ to $10$ biased time-steps with accept/reject stage is a good trade-off between reducing the effect of burn-in in the biased simulation, and efficiency of the resulting mM-MCMC method. We leave an investigation of the bias and the optimal parameters of indirect reconstruction to further research.

\section{Numerical pre-computations using indirect reconstruction} \label{sec:practical}
There are two places in the mM-MCMC algorithm where pre-computations may be needed. First, one may need to pre-compute an approximation to the drift and diffusion terms in the effective dynamics~\eqref{eq:effdyn_intro} to generate the reaction coordinate proposals at the macroscopic level. Similarly, one may want to compute an approximation to the free energy to determine a macroscopic invariant distribution $\bar{\mu}_0$ that approximates the exact time-invariant probability distribution of the reaction coordinates $\mu_0$~\eqref{eq:freeenergy}. 
Second, we need to pre-compute the normalization constant of the indirect reconstruction distribution $N_\lambda$~\eqref{eq:normalizationconstant} for a proper evaluation of  the microscopic acceptance rate~\eqref{eq:matchingacceptance}. Before we can compute this normalization constant, we first need a good approximation of the free energy. Once we have a good approximation to the free energy, this normalization constant is easy to compute.

In Section~\ref{subsec:numeffdyn}, we lay out the computations of the effective dynamics coefficients and the free energy where we propose to use the indirect reconstruction scheme to approximately sample on $\Sigma(z)$. The computations for the normalization constant $N_\lambda$ are explained in Section~\ref{subsec:numfreeenergy}. 

\subsection{Computation of effective dynamics and free energy} \label{subsec:numeffdyn}
Recall from the introduction that the drift and diffusion coefficients $b(z)$ and $\sigma(z)$ are given by the expectations
\begin{equation}\label{eq:eff_coeff}
\begin{aligned}
b(z) &= \mathbb{E}_{\nu(x|z)}[-\nabla V(x) \cdot \nabla \xi(x) + \beta^{-1} \triangle \xi(x)] \\
\sigma(z) &= \mathbb{E}_{\nu(x|z)}[ \norm{\nabla \xi(x)}^2 ], 
\end{aligned}
\end{equation}
computed with respect to the (time-invariant) direct reconstruction distribution $\nu(\cdot| \ z)$~\eqref{eq:nu} on $\Sigma(z)$. Computing these expectations is hard for non-trivial reaction coordinates due to the unknown geometry of $\Sigma(z)$. Existing methods like the projection scheme~\cite{stoltz2010free} or Hamiltonian Monte Carlo~\cite{lelievre2018hybrid} rely on projecting back a proposal to the manifold $\Sigma(z)$, which is achieved by solving a non-linear system for Lagrange multipliers corresponding to the constraint. Solving these non-linear systems can be slow in practice.

As is clear from the companion paper~\cite{vandecasteele2020direct} and Theorem \ref{thm:convergence_matching}, the mM-MCMC method (both with direct and indirect reconstruction) will sample the microscopic distribution $\mu$ exactly, regardless of the choice of macroscopic proposal distribution $q_0$ or (approximate) reaction coordinate distribution $\bar{\mu}_0$. We are therefore willing to tolerate some inaccuracy in the coefficients~\eqref{eq:eff_coeff} to limit the computational overhead on their pre-computation. We propose to sample from the indirect reconstruction distribution $\nu_\lambda(x; z)$, defined in~\eqref{eq:constrainedinvariant}, on the whole configuration space $\mathbb{R}^d$, instead of sampling the direct reconstruction distribution $\nu(x|z)$, defined in~\eqref{eq:nu}, on $\Sigma(z)$. We thus approximate the exact coefficients $b$ and $\sigma$ by
\begin{equation}\label{eq:approx_eff_coeff}
\begin{aligned}
\hat{b}(z) &= \mathbb{E}_{\nu_\lambda(x;z)}[-\nabla V(x) \cdot \nabla \xi(x) + \beta^{-1} \triangle \xi(x)] \\
\hat{\sigma}(z) &= \mathbb{E}_{\nu_\lambda(x;z)}[ \norm{\nabla \xi(x)}^2 ].
\end{aligned}
\end{equation}
Since, in a weak sense, the distribution $\nu_\lambda(\cdot;z)$ lies close to $\nu(\cdot|z)$, the approximate coefficients $\hat{b}$ and $\hat{\sigma}$ will also be a close approximation to the exact coefficients $b$ and $\sigma$, and this approximation is better when $\lambda$ is larger.

For practical computations, we first define a grid of reaction coordinate values $\{z_j\}_{j=1}^J$ on which we approximate the coefficients of the effective dynamics. If the value of $b$ or $\sigma$ is required in an intermediate value $z$, we use linear interpolation between the adjacent grid points to estimate the value at point $z$.  
Starting from a random initial condition, sampling the indirect reconstruction distribution $\nu_\lambda(x;  z_j)$ via a standard MCMC method results in a Markov chain $\{X_j^i\}_{i=1}^N$ whose reaction coordinates lie close to $z_j$. With these MCMC samples, we can then estimate the drift and diffusion terms~\eqref{eq:approx_eff_coeff} as
\begin{equation*}
\begin{aligned}
\hat{b}(z_j) &\approx \frac{1}{N} \sum_{i=1}^N \ -\nabla V(X_j^i) \cdot \nabla \xi(X_j^i) + \beta^{-1} \triangle \xi(X_j^i) \\
\hat{\sigma}^2(z_j) &\approx \frac{1}{N} \sum_{i=1}^N \ \norm{\nabla \xi(X_j^i)}^2.
\end{aligned}
\end{equation*}
The Monte Carlo approximation induces both a statistical error and bias. The statistical error is due to the finite number of samples, and we can make this bias as small as necessary by increasing $N$. There are also two factors that induce a bias with respect to the exact functions $b$ and $\sigma$. The first contribution to the bias is that the microscopic samples $\{X_j^i\}_{i=1}^N$ do not exactly have the prescribed reaction coordinate value $z_j$, since they form samples from $\nu_\lambda(\cdot;z_j)$ instead of $\nu(\cdot|z_j)$. This contribution to the bias can be made as small as desired by increasing the value of $\lambda$. The second bias contribution stems from the fact that we only approximate the effective dynamics on a fixed grid of $z-$values. We can also make this bias term small by defining enough grid points. These errors on the effective dynamics coefficients do not cause a systematic error on the mM-MCMC method, as long as we use the same approximate coefficients in the formulation of the transition probability $q_0$. 

The macroscopic transition distribution, based on the effective dynamics, then reads
\[
q_0(z'|  z) = \left(4 \pi \beta^{-1} \delta t \hat{\sigma}^2(z)\right)^{-n/2} \ \exp\left(-\beta \frac{\norm{z' - z - \hat{b}(z) \delta t}^2}{4 \delta t \hat{\sigma}^2(z)}\right),
\]
where, when $z$ is not a grid point, we define $\hat{b}(z)$ and $\hat{\sigma}(z)$ by linear interpolation between the neighbouring grid cells.

Similarly, we can use the biased stochastic process~\eqref{eq:constrainedsimulation} to compute an approximation to the free energy~\eqref{eq:freeenergy}. We can write the free energy as an expectation over the Lebesgue measure $\sigma_{\Sigma(z)}$
\begin{equation} \label{eq:freeenergyexp}
A(z) =  -\frac{1}{\beta} \ln \left(\mathbb{E}_{\sigma_{\Sigma(z)}(x)}\left[\exp(-\beta V(x)) \norm{\nabla \xi(x)}^{-1}\right]\right).
\end{equation}
We can compute an approximation to the free energy by sampling from the Lebesgue measure $d\sigma_{\Sigma(z_j)}$ for the same grid $\{z_j\}_{j=1}^J$ of reaction coordinate values. We achieve this by again running the biased stochastic process with invariant distribution
\[
\exp\left(-\beta V(x)\right) \exp\left(-\frac{\beta \lambda}{2} (\xi(x) -z)^2\right),
\]
to obtain a set of microscopic samples $\{X_j^i\}_{i=1}^N$. Since the biased simulation scheme actually samples close to $\exp\left(-\beta V(x)\right) \norm{\nabla \xi(x)}^{-1}$ near the sub-manifold $\Sigma(z)$, we associate a normalized weight $w_j^i \propto \exp\left(\beta V(X_j^i)\right) \norm{\nabla \xi(X_j^i)}$ to each microscopic sample $X_j^i$. With these weights, the microscopic samples approximately sample a uniform measure around $\Sigma(z)$ and therefore approximately sample the Lebesgue measure $\sigma_{\Sigma(z)}$ on $\Sigma(z)$ as well.

With the weighted microscopic ensemble $\{X_j^i, w_j^i\}_{i=1}^N$, we subsequently approximate the free energy~\eqref{eq:freeenergyexp} via
\[
A(z_j) \approx \widehat{A}(z_j) = -\frac{1}{\beta} \ln\left( \sum_{i=1}^N w_j^i \exp(-\beta V(X_j^i)) \norm{\nabla \xi(X_j^i)}^{-1} \right).
\]
We can thus use the probability distribution $\bar{\mu}_0 \propto \exp(-\widehat{A}(z))$ as the approximate invariant distribution of the reaction coordinates at the macroscopic level in the mM-MCMC algorithm since we already need an approximation of the free energy to estimate $N_\lambda$.

\subsection{Computing the normalization constant $N_\lambda(\cdot)$} \label{subsec:numfreeenergy}
To close this section, we give a simple numerical scheme to estimate the normalization constant $N_\lambda(\cdot)$ in equation~\eqref{eq:normalizationconstant} of the indirect reconstruction distribution $\nu_\lambda(\cdot ; \cdot)$. This normalization constant can be written as an expected value
\begin{equation} \label{eq:normalizationexp}
N_\lambda(z) = \frac{Z_V}{Z_A} \left(\frac{\lambda \beta}{2\pi}\right)^{n/2}  \mathbb{E}_{\mathcal{N}\left(z, (\lambda \beta)^{-1}\right)} \left[\exp(-\beta A(\cdot))\right],
\end{equation}
but we do not need to compute the normalization constants $Z_V$ and $Z_A$ as these cancel out in the microscopic acceptance rate~\eqref{eq:matchingacceptance}.

We estimate the expected values~\eqref{eq:normalizationexp} at a grid of reaction coordinate values $\{z_j\}_{j=1}^J$ using a Monte Carlo approximation. For every $j$ between $1$ and $J$, we sample $N$ microscopic samples from the $n-$dimensional Gaussian distribution with mean $z_j$ and variance $1/\lambda \beta$. Call these Gaussian samples $\{X_j^i\}_{i=1}^N$. The Monte Carlo approximation of~\eqref{eq:normalizationexp} then reads
\[
N_\lambda(z_j) \approx \frac{Z_V}{Z_A} \left(\frac{\lambda \beta}{2\pi}\right)^{n/2} \frac{1}{N} \sum_{i=1}^N \ \exp\left(-\beta \hat{A}(X_j^i)\right).
\]

We finally note that sampling from the Gaussian distribution $\mathcal{N}(z_j, (\lambda \beta)^{-1})$ is much cheaper than simulating the biased stochastic process to pre-compute the effective dynamics coefficients. It is thus feasible to sample more Gaussian particles to compute the denominator $N_\lambda$ than for the free energy and effective dynamics computations. More Gaussian particles may be necessary since the statistical error on $N_\lambda$ induces a systematic error on the microscopic acceptance probability~\eqref{eq:matchingacceptance}, which we want to be negligible.

\section{Numerical results} \label{sec:results}
In this Section, we numerically investigate the efficiency of the mM-MCMC scheme with indirect reconstruction on two molecular problems: a three-atom molecule and the standard test case alanine-dipeptide. We compare the efficiency gain over the microscopic MALA (Metropolis-adjusted Langevin) method, where we specifically study the impact of the parameters in the indirect reconstruction scheme on this gain. The efficiency gain criterion for a proper comparison of mM-MCMC with the MALA method is explained in Section~\ref{subsec:effcriterion}, and the numerical results for the three-atom molecule and alanine-dipeptide are shown in Sections~\ref{subsec:threeatom} and~\ref{subsec:alanine-dipeptide} respectively.

\subsection{Efficiency criterion} \label{subsec:effcriterion}
Consider a scalar function $F: \mathbb{R}^d \to \mathbb{R}$ and suppose we are interested in the average of $F$ with respect to the Gibbs measure $\mu$,
\[
\mathbb{E}_{\mu}[F] = \int_{\mathbb{R}^d} F(x) \ d\mu(x).
\]
If we sample the invariant measure $\mu$ using an MCMC method, we can estimate the above value as $\hat{F} = N^{-1} \sum_{n=1}^N F(x_n)$ with an ensemble of microscopic samples $\{x_n\}_{n=1}^N$. The variance on this estimate is
\begin{equation} \label{eq:varF}
\text{Var}[\hat{F}] = \frac{\sigma_F^2 \ K_{\text{corr}}}{N},
\end{equation}
where $\sigma_F^2$ is the inherent variance of $F$, 
\[
\sigma_F^2 = \int_{\mathbb{R}^d} \ \left(F(x) - \mathbb{E}_{\mu}[F] \right)^2 d\mu(x),
\]
and $K_{\text{corr}}$ is defined as
\[
K_{\text{corr}} = 1 + \frac{2}{\sigma_F^2} \ \sum_{n=1}^N \ \mathbb{E}\left[(F(x_n) - \mathbb{E}_{\mu}[F(x_n)])(F(x_0) - \mathbb{E}_{\mu}[F(x_0)] )\right],
\]
with $x_0$ the initial value of the Markov chain~\cite{meyn2012markov} .

The extra factor $K_{\text{corr}}$ in~\eqref{eq:varF} arises because the Markov chain Monte Carlo samples are not independent. The higher $K_{\text{corr}}$, the more correlated the MCMC samples and the higher the variance~\eqref{eq:varF}. Another interpretation of the correlation parameter $K_{\text{corr}}$ is that the `effective' number of samples is $N/ K_{\text{corr}}$.

In the following numerical experiments, we are interested in reducing the variance on the estimator $\hat{F}$ with mM-MCMC using the same runtime, compared the microscopic MALA algorithm. Equivalently, we want to increase the effective number of samples $N/K_{\text{corr}}$ for a fixed runtime. We therefore define the efficiency gain of mM-MCMC over the microscopic MALA method as
\begin{equation} \label{eq:effgain}
\text{Gain} = \frac{\text{Var}[\hat{F}]_{\text{micro}}}{\text{Var}[\hat{F}]_{\text{mM}}} \ \frac{T_{\text{micro}}^N}{T_{\text{mM}}^N} = \frac{K_{\text{corr, micro}}}{K_{\text{corr, mM}}} \ \frac{T_{\text{micro}}^N}{T_{\text{mM}}^N}.
\end{equation}
Here, $T_{\text{micro}}^N$ is the measured execution time of the microscopic MCMC method for a fixed number of sampling steps $N$ and $T_{\text{mM}}^N$ is the measured execution time for the same number of steps $N$ of the mM-MCMC scheme. 

Usually, the execution time $T_{mM}^N$ of the mM-MCMC method with indirect reconstruction is larger than the execution time $T_{micro}^N$ for the MALA method for the same number of sampling steps $N$, due to the overhead of the indirect reconstruction scheme. Indeed, when a reaction coordinate value is accepted at the macroscopic level, we take $K$ time steps of the biased simulation scheme where we evaluate the potential energy and its gradient at every biased time step. The exact increase in execution time for the same number of sampling steps $N$ depends on the value of the macroscopic acceptance rate. We will show that the effective number of samples $N/K_{\text{corr, mM}}$ of the mM-MCMC scheme is usually orders of magnitude higher than the effective number of samples of the MALA scheme, $N/K_{\text{corr, micro}}$. By first sampling a reaction coordinate value with a large time step $\Delta t$ at the macroscopic level, the correlation between two reaction coordinate values will, on average, be lower than the correlation between two microscopic samples generated with a small time step $\delta t$ by the MALA method. Currently, however, we have no analytic formulas linking the correlations at the macroscopic and microscopic levels, so we will demonstrate this claim numerically. Combining both the decrease in execution time for the same number of sampling steps and the increase of the effective number of samples, we expect that mM-MCMC will be able to gain over the MALA method for moderate to large time-scale separations. We will show that the higher the time-scale separation, the larger the efficiency gain will be.

\subsection{A three-atom molecule} \label{subsec:threeatom}
\begin{figure}
	\centering
	\includegraphics[width=0.3\linewidth]{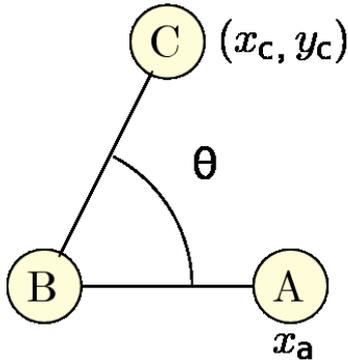}
	\caption{The three-atom molecule. Atom $A$ is constraint to the $x$-axis with $x$-coordinate $x_a$, atom $B$ is fixed at the origin of the plane and atom $C$ lies om the two-dimensional plane with Cartesian coordinates $(x_c, y_c)$.}
	\label{fig:triatommolecule2}
\end{figure}

\paragraph{\textbf{Model problem}}
In this section, we consider the mM-MCMC algorithm on a simple, academic, three-atom molecule, as first introduced in~\cite{legoll2010effective}. The three-atom molecule has a central atom $B$, that we fix at the origin of the two-dimensional plane, and two outer atoms, $A$ and $C$. To fix the superfluous degrees of freedom, we constrain atom $A$ to the $x-$axis, while $C$ can move freely in the plane. The three-atom molecule is depicted on Figure~\ref{fig:triatommolecule2}. 

The potential energy for the three-atom system consists of three terms,
\begin{equation} \label{eq:trheeatompotential}
V(x_a, x_c, y_c) = \frac{1}{2\varepsilon} \ (x_a-1)^2 + \frac{1}{2\varepsilon} \ (r_c-1)^2 + \frac{208}{2}\left(\left(\theta-\frac{\pi}{2}\right)^2 - 0.3838^2\right)^2,
\end{equation}
where $x_a$ is the $x-$coordinate of atom $A$ and $(x_c, y_c)$ are the Cartesian coordinates of atom $C$. The bond length $r_c$ between atoms $B$ and $C$ and the angle $\theta$ between atoms $A$, $B$ and $C$ are defined as
\begin{equation*}
\begin{aligned}
r_c &= \sqrt{x_c^2 + y_c^2} \\ 
\theta &= \arctantwo (y_c, x_c).
\end{aligned}
\end{equation*}
The first term in~\eqref{eq:trheeatompotential} describes the vibrational potential energy of the bond between atoms $A$ and $B$, with equilibrium length $1$. Similarly, the second term describes the vibrational energy of the bond between atoms $B$ and $C$ with bond length $r_c$. Finally, the third term determines the potential energy of the angle $\theta$ between the two outer atoms, which has an interesting bimodal behaviour. The distribution of $\theta$ has equilibrium values, one at $\frac{\pi}{2} -0.3838$ and another at $\frac{\pi}{2} + 0.3838$. 

In the following set of experiments, we define the angle $\theta$ as reaction coordinate, i.e., 
\[
\xi(x) = \theta,
\]
and we study the efficiency of mM-MCMC for this reaction coordinate choice. The angle $\theta$ between the two outer atoms is purely slow since its corresponding term in the potential energy~\eqref{eq:trheeatompotential} is independent from $\varepsilon$. We therefore expect to obtain significant efficiency gain using the mM-MCMC method with this reaction coordinate function.

\paragraph{Overview of this section}
In this section, we perform the following numerical experiments on the three-atom molecule. In Section~\ref{subsubsec:threeatomprecomp}, we carry out the pre-computations of the free energy and the effective dynamics of $\xi$. Subsequently, in Section~\ref{subsubsec:visual_inspection}, we visually inspect the efficiency gain of the mM-MCMC method with indirect reconstruction over the microscopic MALA method. In Section~\ref{subsubsec:triatom_angle}, we numerically investigate the efficiency gain on the estimated mean and variance for multiple values of $\varepsilon$, followed by a comparison between the direct and indirect reconstruction variants in Section~\ref{subsubsec:badrc} for the same values of $\varepsilon$. Finally, we determine, numerically, the optimal value for the strength of the biasing potential, $\lambda$, that maximizes the efficiency gain of mM-MCMC with reaction coordinate $\xi$ for multiple values of $\varepsilon$ in Section~\ref{subsubsec:threeatomparam}.

\subsubsection{Pre-computations of the free energy and effective dynamics} \label{subsubsec:threeatomprecomp}
\paragraph{Experimental setup}
In this section, we numerically approximate the free energy and the coefficients in the effective dynamics of reaction coordinate $\xi_1$ using the indirect reconstruction scheme, as outlined in Section~\ref{sec:practical}. We define $200$ grid points between $0$ and $\pi$. The parameters for the indirect reconstruction scheme are $\lambda = 100/\varepsilon$, $\delta t = 1/\lambda$ and we use $N=10000$ microscopic Monte Carlo samples per grid value of the reaction coordinate. The inverse temperature is fixed to $\beta=1$. 

\begin{figure}
	\centering
	\begin{subfigure}[b]{0.5\textwidth}
		\centering
		\includegraphics[width=1\linewidth]{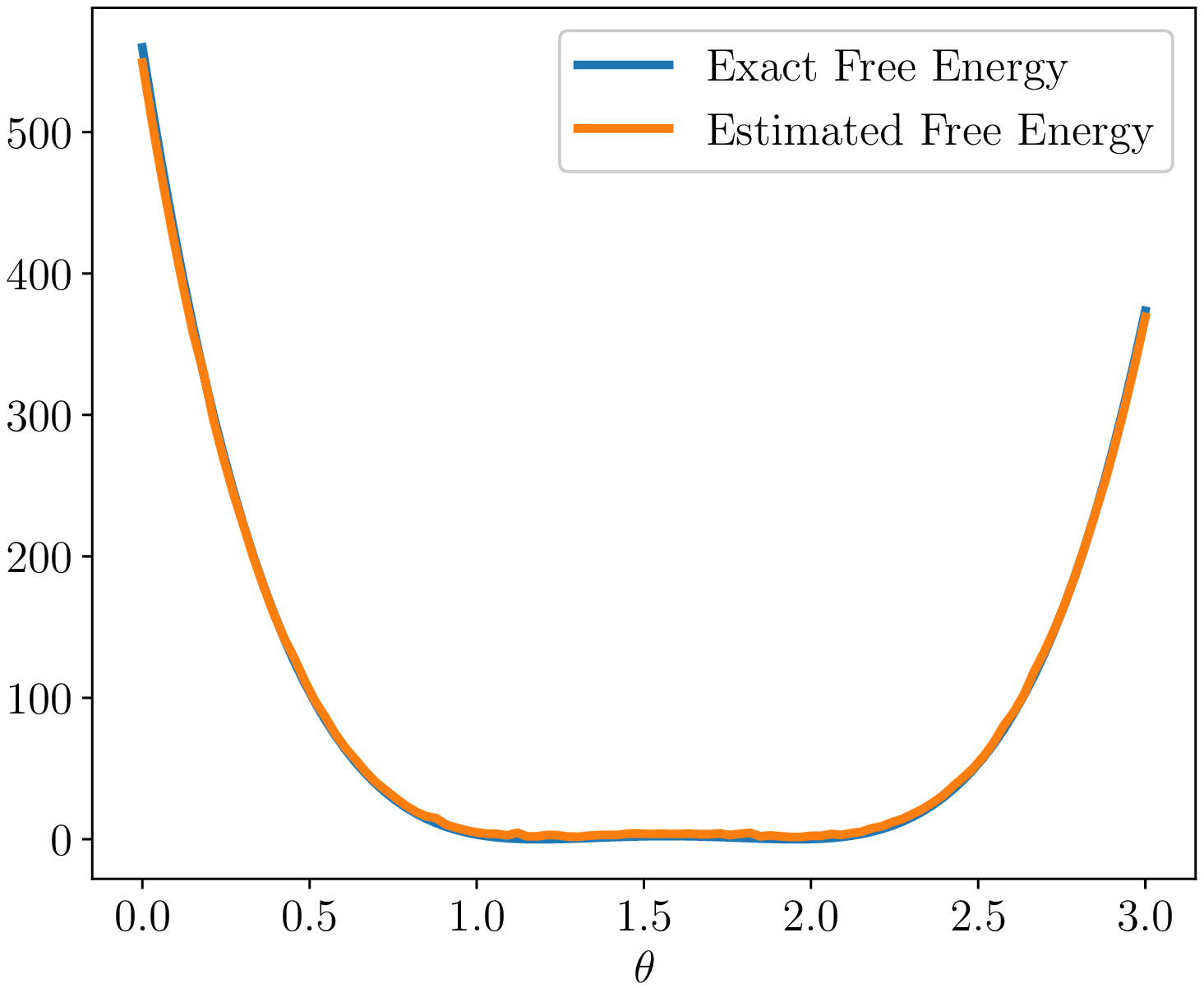}
	\end{subfigure}%
	\begin{subfigure}[b]{0.5\textwidth}
		\centering
		\includegraphics[width=1\linewidth]{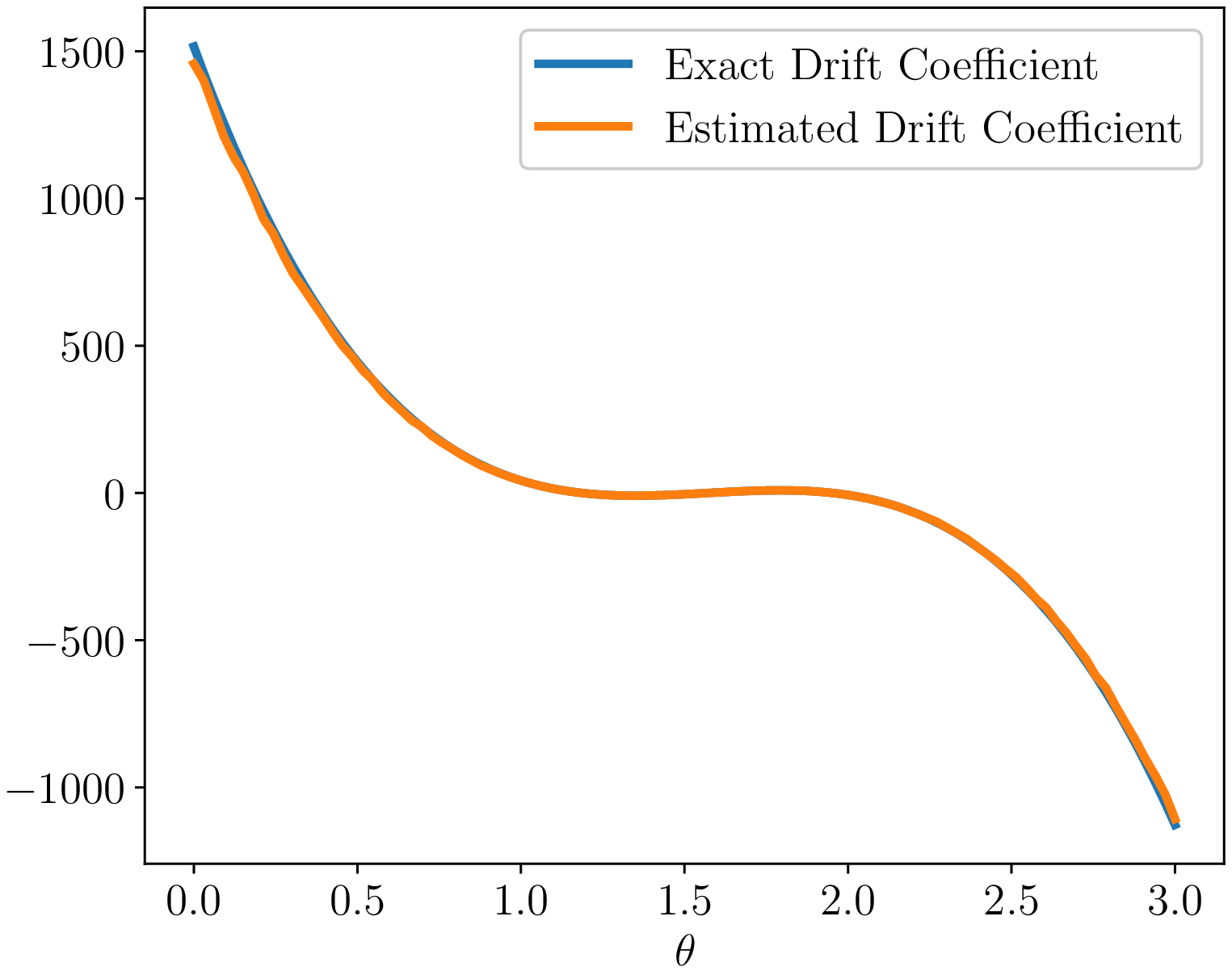}
	\end{subfigure}
	\caption{Estimated free energy (left) and drift coefficient in the effective dynamics (right) using $100$ reaction coordinate grid points and $N=10000$ microscopic samples per reaction coordinate value. The numerical solutions are compared with the analytical expressions~\eqref{eq:freedriftanalytic} (blue). Based on the numerical results, the numerical scheme based on indirect reconstruction proposed in Section~\ref{sec:practical} yields a good approximation to the free energy and effective dynamics.}
	\label{fig:free_energy_drift}
\end{figure}

In case of reaction coordinate $\xi_1$, the free energy $A(\theta)$ and the drift coefficient $b(\theta)$ are analytically known, i.e.,
\begin{equation} \label{eq:freedriftanalytic}
\begin{aligned}
A(\theta) &= \frac{k}{2}\left(\left(\theta-\frac{\pi}{2}\right)^2 -\delta \theta^2\right)^2 \\
b(\theta) &= -2 k \left(\left(\theta-\frac{\pi}{2}\right)^2 -\delta \theta^2\right) \left(\theta-\frac{\pi}{2}\right),
\end{aligned}
\end{equation}
and the diffusion coefficient $\sigma(\theta)$ is always $1$ because $\norm{\nabla \xi_1}=1$. In Figure~\ref{fig:free_energy_drift}, we numerically compare the estimated free energy and drift coefficient to the analytic formulas~\eqref{eq:freedriftanalytic}.

\paragraph{Numerical results}
Theoretically, there are two error contributions to the pre-computations: a statistical error due to the finite number of microscopic samples $N$ to approximate $A(\theta)$ and $b(\theta)$, and a deterministic bias because $\lambda$ is finite. As one can see from the numerical results, the estimated free energy and estimated drift coefficient lie close to their corresponding analytic values.  Since there are almost no fluctuations on both curves, the dominant error contributor appears to be the bias due to $\lambda$. Indeed, the farther the value of $\theta$ lies from $\pi/2$, the larger the absolute error becomes, although relatively, both values lie close to each other. We will use these estimated quantities for the remaining numerical experiments with reaction coordinate $\xi_1$ in this section.

\subsubsection{Visual inspection of mM-MCMC with indirect reconstruction} \label{subsubsec:visual_inspection}
\paragraph{Experimental setup}
Now that we have an accurate approximation to the free energy and the effective dynamics coefficients, we can visually inspect how the mM-MCMC algorithm with indirect reconstruction compares to the microscopic MALA method for a fixed value of $\varepsilon=10^{-6}$. When $\varepsilon$ is small, the microscopic time step $\delta t$ of the microscopic MALA method must be small as well since $\delta t$ scales with $\varepsilon$ due to the stiffness of the problem. As a result, we expect that the microscopic MALA method will remain stuck for a long time in one of the potential wells of the macroscopic variable $\theta$.  

\begin{figure}
	\centering
	\begin{subfigure}[b]{0.5\textwidth}
		\centering
		\includegraphics[width=0.9\linewidth]{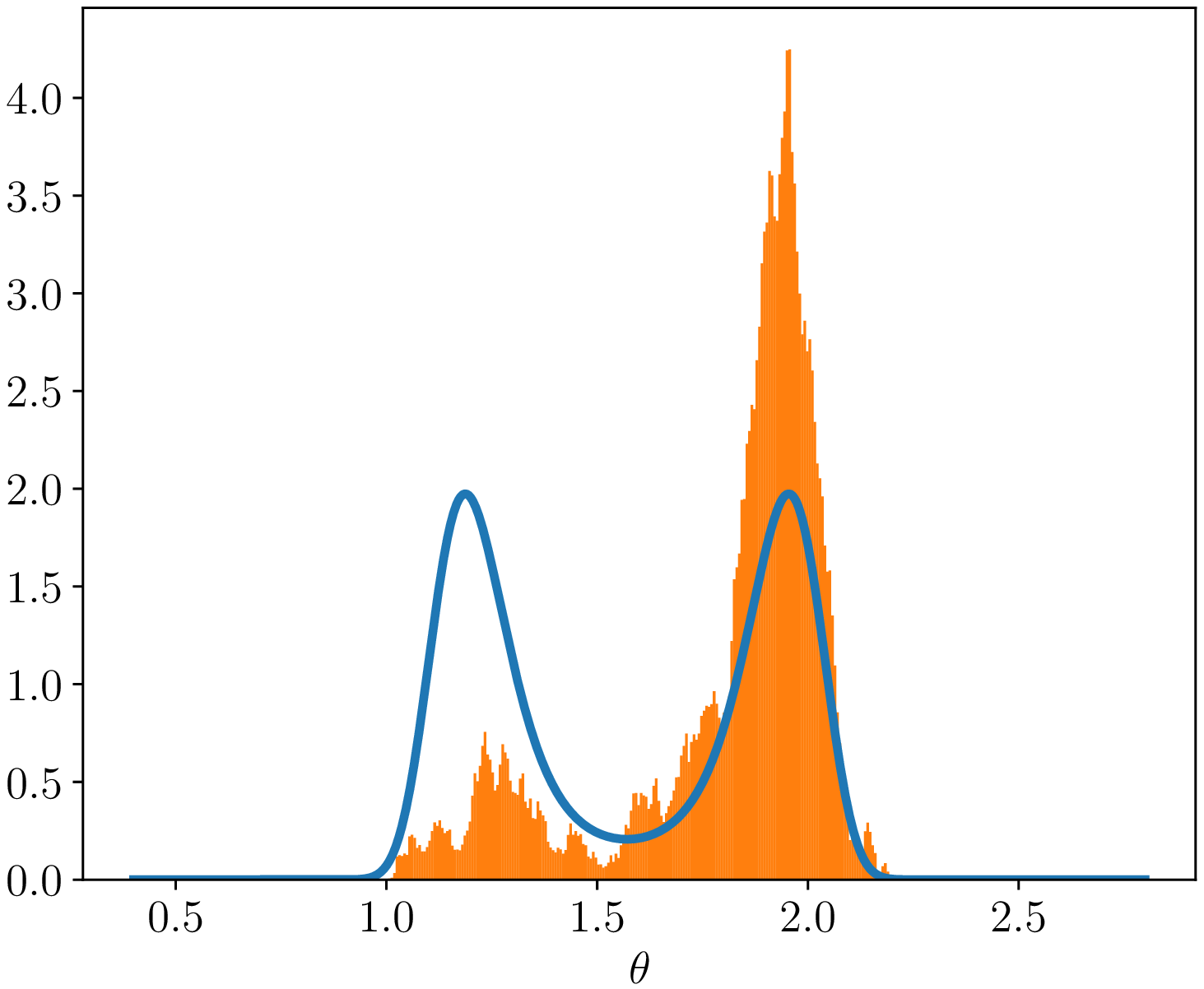}
	\end{subfigure}%
	\begin{subfigure}[b]{0.5\textwidth}
		\centering
		\includegraphics[width=0.9\linewidth]{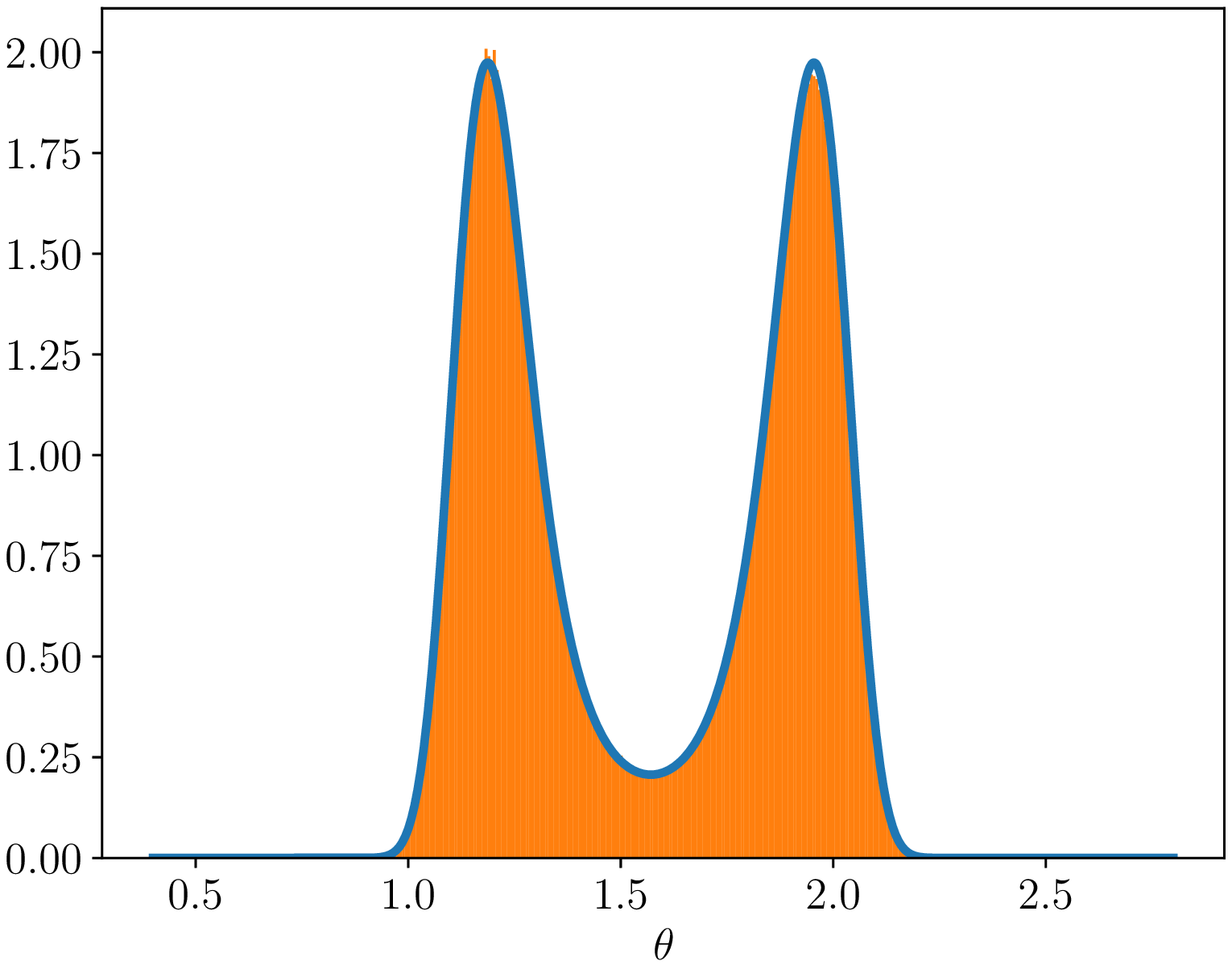}
	\end{subfigure}
	\caption{Visual representation of the histogram of $\theta$ of the microscopic MCMC method (left) and mM-MCMC (right) on the three-atom molecule, with reaction coordinate $\theta$. The simulation parameters are $\varepsilon=10^{-6}$ and the number of samples is $N=10^6$. The microscopic MCMC method remains stuck in the potential well of $\theta$ around $1.95$, while the mM-MCMC method samples the distribution well.}
	\label{fig:threeatom_mM_micro}
\end{figure}

In this experiment, we base the macroscopic proposals on the Euler-Maruyama discretization of the approximate effective dynamics from the previous section with step size $\Delta t = 0.01$. The macroscopic invariant distribution $\bar{\mu}_0$ is also based on the approximate free energy that we computed in Section~\ref{subsubsec:threeatomprecomp}. The indirect reconstruction parameters are $K=5$, $\lambda = \varepsilon^{-1}$ and the biased step size is $\delta t=\varepsilon$. This value is also the step size of the microscopic MALA method. The numerical results are shown in Figure~\ref{fig:threeatom_mM_micro} for the microscopic MALA method (left) and mM-MCMC (right).

\paragraph{Numerical results}
As one can see on Figure~\ref{fig:threeatom_mM_micro}, the microscopic MALA method over-samples the right potential well of $\theta$ and puts only a small fraction of the microscopic samples in the left potential well due to the small time steps $\delta t = \varepsilon$. However, the mM-MCMC method with indirect reconstruction is able to sample the distribution of $\theta$ accurately due to the large macroscopic time steps. Therefore, we may expect a large efficiency gain of mM-MCMC over the microscopic MALA method when estimating some quantities of interest of $\theta$.

\subsubsection{Efficiency gain as a function of time-scale separation} \label{subsubsec:triatom_angle}

\begin{table}
	\centering
	\begin{tabular}{c|c|c|c|c|c}
		\centering
		$\varepsilon$ & \pbox{15cm}{Macroscopic \\ acceptance rate} & \pbox{15cm}{Microscopic \\ acceptance rate} & \pbox{15cm}{Runtime \\ gain} & \pbox{15cm}{Variance \\ gain} & \pbox{15cm}{Total \\ efficiency gain} \\
		\hline
		$10^{-3}$ & 0.749975           &            0.993528      &     0.212558       &      10.3927         &         2.20905\\
		$10^{-4}$ & 0.74936          &              0.993588     &      0.211415      &       68.6401          &       14.5115 \\
		$10^{-5}$ & 0.750197          &             0.993299    &       0.212561     &       920.651        &         195.695 \\
		$10^{-6}$ & 0.749498              &         0.993405  &         0.237227   &        7041.69          &       1670.48 \\
	\end{tabular}
	\caption{A summary of different statistics of the mM-MCMC method when applied to estimating the mean of $\theta$, for multiple values of $\varepsilon$. First of all, note that the macroscopic acceptance rate remains constant when $\varepsilon$ decreases, which is a logical consequence from the fact that $\theta$ is independent of the small scales. Second, the microscopic acceptance rate is very close to $1$, showing that the indirect reconstruction step is close to the direct reconstruction distribution~\eqref{eq:nu}. Additionally, a high microscopic acceptance rate indicates that almost no reconstructed microscopic samples are rejected, making the indirect reconstruction scheme efficient. Third, note that the execution time of mM-MCMC is higher than the execution time of the microscopic MALA method due to the computational overhead of the biased simulation. Fortunately, the increase in execution time is independent of $\varepsilon$. However, the decrease in variance of the estimate mean of $\theta$ increases linearly with decreasing $\varepsilon$ and so does the efficiency gain.}
	\label{tab:threeatomgain_mean}
\end{table}

\begin{figure}
	\centering
	\includegraphics[width=0.55\linewidth]{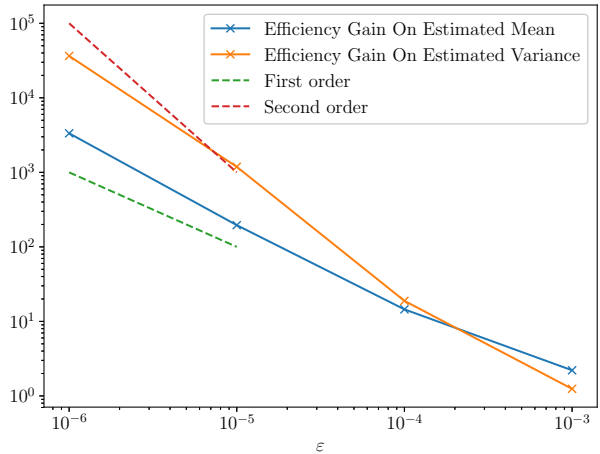}
	\caption{Efficiency gain of mM-MCMC over standard MCMC as a function of the time-scale separation. The efficiency gain for the mean of $\theta$ (blue) increases linearly with decreasing $\varepsilon$, while the efficiency gain for the variance of $\theta$ increases faster than linearly, although slower than quadratically.}
	\label{fig:gainepstriatom}
\end{figure}

\paragraph{Experimental setup}
Following the previous experiment, we numerically compare the efficiency gain of the mM-MCMC method over the microscopic MALA method for different values of $\varepsilon$. We let the small-scale parameter $\varepsilon$ vary between $10^{-6}$ and $10^{-3}$ and we consider the efficiency gain on the estimated mean and the estimated variance of $\theta$. To measure the efficiency gain, we use the efficiency gain criterion that we explained in Section~\ref{subsec:effcriterion}, where we average the estimated mean and the estimated variance over $100$ independent runs for a good assessment of the efficiency gain. The timings were performed on a Xeon Gold 6140 (Skylake) processor and the numerical parameters are the same as in the previous experiment. The efficiency gains on the estimated mean and variance of $\theta$ are depicted on Figure~\ref{fig:gainepstriatom} as a function of $\varepsilon$. We also display the averaged macroscopic and microscopic acceptance rates, the reduction in runtime, as well as the reduction of the variance on the estimated quantities of interest obtained by mM-MCMC with indirect reconstruction in Table~\ref{tab:threeatomgain_mean} where we consider the estimated mean of $\theta$, and in Table~\ref{tab:threeatomgain_variance} for the estimated variance of $\theta$.

\paragraph{Numerical results}
First, on Figure~\ref{fig:gainepstriatom}, we see that the efficiency gain of mM-MCMC increases at least linearly with decreasing $\varepsilon$. For large values of $\varepsilon$, one can see that there is almost no efficiency gain at all. The reason for this behaviour is clear from Tables~\ref{tab:threeatomgain_mean} and~\ref{tab:threeatomgain_variance}. From the third column of these tables, we conclude that the runtime of mM-MCMC with indirect reconstruction is larger than that of the microscopic MALA method due to the computational overhead of the biased simulation. Hence, for large values of $\varepsilon$, this computational overhead (third column) is dominant over the reduction in variance (fourth column). However, when $\varepsilon$ is small, the reduction in variance by mM-MCMC is dominant over the increase in computational time. Currently, we do not have analytic expressions for the gain as a function of the time-scale separation, but Figure~\ref{fig:gainepstriatom} clearly shows the merit of mM-MCMC with indirect reconstruction for medium and large time-scale separations.

\begin{table}
	\centering
	\begin{tabular}{c|c|c|c|c|c}
		\centering
		$\varepsilon$ & \pbox{15cm}{Macroscopic \\ acceptance rate} & \pbox{15cm}{Microscopic \\ acceptance rate} & \pbox{15cm}{Runtime \\ gain} & \pbox{15cm}{Variance \\ gain} & \pbox{15cm}{Total \\ efficiency gain} \\
		\hline
		$10^{-3}$ & 0.749975              &         0.993528  &         0.212558     &        5.85671         &         1.24489\\
		$10^{-4}$ & 0.74936     &                  0.993588      &     0.211415     &       89.2989                &  18.8791 \\
		$10^{-5}$ & 0.750197         &              0.993299        &   0.212561    &      5580.75      &            1186.25 \\
		$10^{-6}$ & 0.749498          &             0.993405      &     0.237227  &      153706           &         36463.2 \\
	\end{tabular}
	\caption{A summary of different statistics of the mM-MCMC method when applied to estimating the variance of $\theta$, for multiple values of $\varepsilon$. The conclusions on the macroscopic and microscopic acceptance rates and the increase in execution time are the same as in Table~\ref{tab:threeatomgain_mean}.}
	\label{tab:threeatomgain_variance}
\end{table}

As a second observation, note that the macroscopic acceptance rates (first column of both tables) is independent of $\varepsilon$, as we intuitively may expect. Also, we see that the microscopic acceptance rate after indirect reconstruction (second column) is close to $1$ so that only little redundant computational work is performed during indirect reconstruction. This result shows that the indirect reconstruction is an efficient technique to reconstruct a microscopic sample close to a sub-manifold of constant reaction coordinate value.

\subsubsection{Comparison of mM-MCMC with direct and indirect reconstruction}\label{subsubsec:badrc}

\begin{figure}
	\centering
	\includegraphics[width=0.55\linewidth]{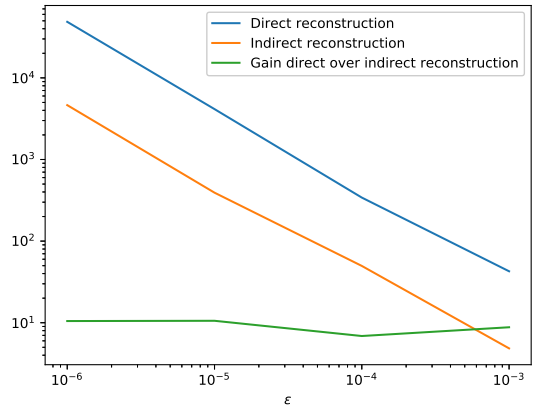}
	\caption{Efficiency gain of mM-MCMC with direct reconstruction (blue) and indirect reconstruction (orange) over the microscopic MALA method as a function of $\varepsilon$. The green curve measures the efficiency gain of the direct reconstruction algorithm over the indirect reconstruction method. One can see that the efficiency gain is a constant factor of less than $10$ lower then the gain made by the direct reconstruction variant.}
	\label{fig:directindirect}
\end{figure}

\paragraph{Experimental setup}
Having studied the performance of mM-MCMC with indirect reconstruction for multiple values of the time-scale separation, we now compare the performance of this method to its direct reconstruction variant. We expect the direct reconstruction algorithm to be faster for a given number of sampling steps due to the computational overhead of the biased dynamics, while the reduction in variance on estimated quantities should be almost identical. In Figure~\ref{fig:directindirect}, we depict the efficiency gain of both mM-MCMC variants over the microscopic MALA method on the estimated mean of $\theta$ for a large range of values of $\varepsilon$. In addition, we also plot the efficiency gain of the direct reconstruction algorithm over \emph{the indirect reconstruction algorithm} for $N=10^6$ sampling steps. We fix the macroscopic time step at $\Delta t = 0.02$ for both mM-MCMC variants for a good comparison. The numerical parameters for the indirect reconstruction method are $K = 5$, $\lambda = \varepsilon^{-1}$ and $\delta t = \varepsilon$, and the time step of the microscopic MALA method is also $\varepsilon$. In Table~\ref{tab:directindirect}, we also show the gain in runtime, the gain in variance on the estimated mean of $\theta$ and the total efficiency gain of mM-MCMC with direct reconstruction over mM-MCMC with indirect reconstruction.

\paragraph{Numerical results}
The efficiency gain of mM-MCMC with indirect reconstruction is a constant factor lower than the efficiency gain of its direct reconstruction variant, independent of the time-scale separation. If we diagnose this effect more carefully in Table~\ref{tab:directindirect}, one can see that the lower efficiency gain is purely due to the larger runtime of mM-MCMC with indirect reconstruction. The decrease in variance on the estimated mean of $\theta$ is almost the same.

\begin{table}
	\centering
	\begin{tabular}{c|c|c|c}
		\centering
		$\varepsilon$  & \pbox{15cm}{Runtime \\ gain} & \pbox{15cm}{Variance \\ gain} & \pbox{15cm}{Total \\ efficiency gain} \\
		\hline
		$10^{-3}$  &         10.0329     &     0.877206        &    8.80089\\
		$10^{-4}$  &   10.1001     &     0.680705     &       6.87522 \\
		$10^{-5}$  &   10.1374     &     1.04204   &         10.5635 \\
		$10^{-6}$  &    9.99052 &        1.57429     &       10.4853 \\
	\end{tabular}
	\caption{A summary of different statistics that summarize the efficiency gain of mM-MCMC with direct reconstruction over the indirect reconstruction algorithm. The efficiency gain of the direct reconstruction algorithm is almost the same for a large range of values of $\varepsilon$ and this efficiency gain is completely due to the lower runtime of the direct reconstruction. Both variants obtain the exact same variance reduction over the microscopic MALA method.}
	\label{tab:directindirect}
\end{table}

\subsubsection{Impact of $\lambda$ on the efficiency of mM-MCMC} \label{subsubsec:threeatomparam}
\begin{figure}
	\centering
	\includegraphics[width=0.85\linewidth]{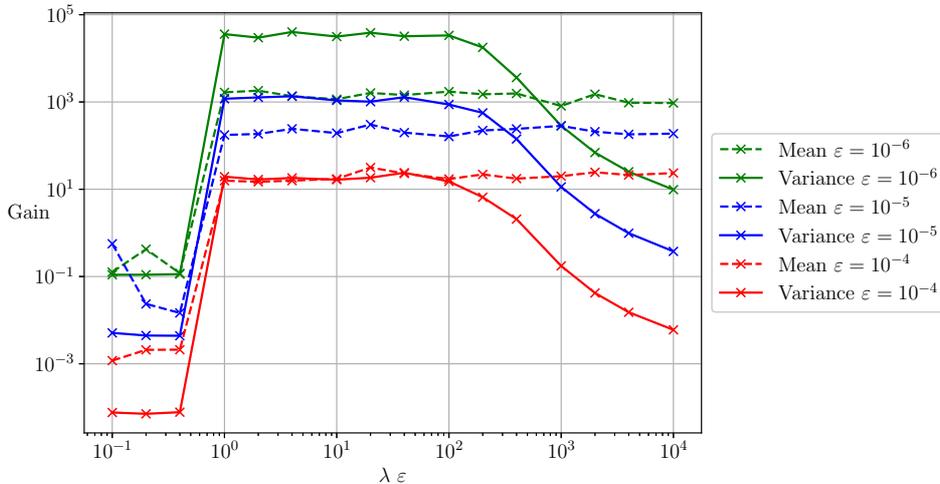}
	\caption{Efficiency gain of mM-MCMC on the estimate of the mean (full lines) and variance (dotted lines) of $\theta$, for different values of the time-scale separation. When $\lambda$ is larger $1/\varepsilon$ there is a clear efficiency gain. This gain is almost constant for a large range of values for $\lambda$. When $\lambda$ is too large, we need more biased steps to come near to the sampled value at the macroscopic level, and when $\lambda$ is smaller than $1/\varepsilon$, we take the sampled reaction value not enough into account and loose efficiency.}
	\label{fig:efflambda}
\end{figure}

\paragraph{Experimental setup}
For the final experiment on the three-atom molecule, we investigate the effect of the magnitude of $\lambda$ on the efficiency gain of the estimated mean and the estimated variance of $\theta$. Specifically, we compute the efficiency gain of mM-MCMC with indirect reconstruction over the microscopic MALA method for three values of the time-scale separation: $\varepsilon = 10^{-4}, \ 10^{-5}, \ \text{and} \  10^{-6}$. For a good comparison, we keep the number of biased time steps fixed at $K = 5$ and the step size of the biased simulation fixed at $\delta t=\lambda^{-1}$. The other numerical parameters are $\beta=1$, the macroscopic time step is $\Delta t = 0.01$, the time step of the microscopic MALA method is $\varepsilon$ and we take $N=10^6$ sampling steps. Moreover, we compute the efficiency gain by averaging the estimated quantities over 100 independent runs. On Figure~\ref{fig:efflambda}, we show the efficiency gain of mM-MCMC as a function of $\lambda$, for multiple values of $\varepsilon$. Also, in Table~\ref{tab:threeatomgain_lambda}, we gather some statistics of the mM-MCMC method for the estimated variance of $\theta$ for several values of $\lambda$ and a fixed value of $\varepsilon=10^{-6}$.

\begin{table}[!h]
	\centering
	\begin{tabular}{c|c|c|c|c|c}
		\centering
		$\lambda  \cdot \varepsilon$ & \pbox{15cm}{Macroscopic \\ acceptance rate} & \pbox{15cm}{Microscopic \\ acceptance rate} & \pbox{15cm}{Runtime \\ gain} & \pbox{15cm}{Variance \\ gain} & \pbox{15cm}{Total \\ efficiency gain} \\
		\hline
		$0.1$ &                0.749749   &                    0.993677    &       0.218974         &   0.500816          &       0.109666\\
		$1$ & 0. 0.74966                     &   0.993468     &      0.233651  &     153735        &            35920.4 \\
		$10$ & 0.7496                         &0.993463      &     0.22258       & 141816            &        31565.4 \\
		$100$ & 0.749546                     &  0.993421     &      0.219963  &     152517           &         33548.1 \\
		$1000$ & 0.749553        &               0.993469    &       0.218501  &       1303.67       &            284.853 \\
	\end{tabular}
	\caption{Several statistics of the performance of mM-MCMC with indirect reconstruction over the microscopic MALA method when estimating the variance of $\theta$, for several values of $\lambda$. The time-scale separation parameter is fixed at $\varepsilon=10^{-6}$. The macroscopic and microscopic acceptance rates and the increase in execution time remain constant with varying $\lambda$. As is also visible from the greed solid line on Figure~\ref{fig:efflambda}, when $\lambda < \varepsilon^{-1}$, the efficiency gain is low because there is no gain on the variance at all. However, from the point $\lambda > \varepsilon^{-1}$, there is a significant variance reduction on the estimated variance of $\theta$, and therefore a significant increase in efficiency gain. Finally, when $\lambda$ is too large, relative to the fastest modes in the potential energy of the system, the total efficiency gain decreases again due to an increase of the variance on the estimated variance of $\theta$ by mM-MCMC.}
	\label{tab:threeatomgain_lambda}
\end{table}

\paragraph{Numerical results}
First, note that the efficiency gain increases when $\varepsilon$ decreases, an effect that we already studied in Section~\ref{subsubsec:triatom_angle}. Second, when $\lambda < 1/\varepsilon$ there is no efficiency gain at all. Indeed, as we intuitively mentioned in Section~\ref{subsec:optimalparameters}, when $\lambda$ is smaller than the stiffest mode in the system, the reaction coordinate in the biased dynamics will be less driven towards the value sampled at the macroscopic level. We are then effectively ignoring the macroscopic MCMC step so that there is no variance reduction on expectations of $\theta$. This effect is also visible in the first row of Table~\ref{tab:threeatomgain_lambda} where there is indeed no gain on the variance on the estimated variance of $\theta$. On the other hand, when $\lambda$ is very large (larger than $10^2 / \varepsilon$ in this case), the efficiency also starts to decrease as we are not simulating enough biased steps for the reaction coordinate value to approximate the sampled value at the macroscopic level well. This effect is especially visible on the efficiency gain of the estimate variance of $\theta$ in the bottom row of Table~\ref{tab:threeatomgain_lambda}. We would therefore need more than $5$ biased steps to equilibrate around each sampled reaction coordinate value at the macroscopic level of the mM-MCMC algorithm, also reducing the efficiency. To conclude, there is a large range of values for $\lambda$ that give a large and almost identical efficiency gain (middle rows of Table~\ref{tab:threeatomgain_lambda}). This range, between $1/\varepsilon$ and $100/\varepsilon$, is the same for a large range of $\varepsilon$ values, supporting our claim in Section~\ref{subsec:optimalparameters}. In practice, it is therefore a good idea to choose $\lambda$ approximately on the order of the stiffest mode of the molecular system.

\subsection{Alanine-dipeptide} \label{subsec:alanine-dipeptide}
\begin{figure}
	\centering
	\includegraphics[width=0.5\linewidth]{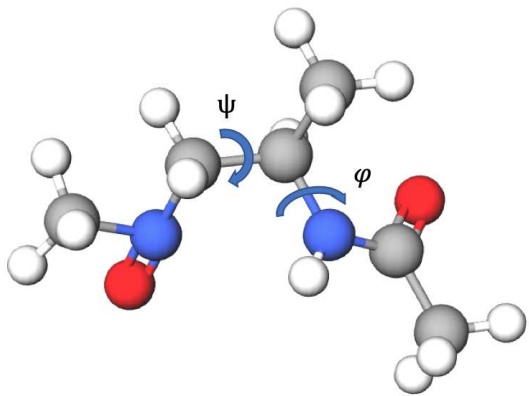}
	\caption{The alanine-dipeptide molecule. The carbon atoms are grey, hydrogen white, nitrogen blue and oxygen is red.}
	\label{fig:alanine-dipeptide}
\end{figure}

\paragraph{Model problem}
For the second molecular test case, we consider the alanine-dipeptide molecule, a much-used test case for numerical methods in molecular dynamics~\cite{leimkuhler2006new}. Alanine-dipeptide has two internal torsion angles, $\phi$ and $\psi$, of which $\psi$ is slow and $\phi$ is an order of magnitude faster than $\psi$. The molecule is depicted on Figure~\ref{fig:alanine-dipeptide}. The first torsion angle $\phi$ determines the position of the right subgroup of the molecule, starting from the carbon-nitrogen bond with the central carbon molecule. The other torsion angle $\psi$ then determines the position of the left subgroup, starting from the carbon-carbon bond with the central carbon atom.

The potential energy function of alanine-dipeptide contains of a term for each atom-atom bond, a term for each angle between two bonds that are connected by a common atom and a term for each of the two torsion angles. We do not consider non-bonded interactions. 
The terms and parameters of the potential energy are summarised in Table~\ref{tab:alaninedipeptide}. Note that each potential energy term corresponding to a torsion angle has only one stable conformation. The parameters in the potential energy are taken from~\cite{head1991strategy}.

\begin{table}
	\centering
	\begin{tabular}{c|c|c}
		Term & form & parameters \\
		\hline
		C-C Bond & $0.5 \ k_{CC} \ (r - r_{CC})^2$ & $k_{CC} = 1.17 \cdot 10^{6}$, $r_{CC} = 1.515$ \\
		C-N Bond & $0.5 \ k_{CN} \ (r - r_{CN})^2$ & $k_{CN} = 1.147 \cdot 10^{6}$, $r_{CN} = 1.335$ \\
		C-C-N Angle & $0.5 \ k_{CCN} \ (\theta - \theta_{CCN})^2$ & $k_{CCN} = 2.68 \cdot 10^{5}, \ \theta_{CCN} = 113.9 \deg$ \\
		C-N-C Angle & $0.5 \ k_{CNC} \ (\theta  - \theta_{CNC})^2$ & $k_{CNC} = 1.84 \cdot 10^5, \ \theta_{CNC} = 117.6 \deg$ \\
		Torsion Angle $\phi$ & $k_{\phi} (1 + \cos(\phi + \pi))$ & $k_{\phi} = 3.98 \cdot 10^4 $ \\
		Torsion Angle $\psi$ & $k_{\psi} (1 + \cos(\psi + \pi))$ & $k_{\psi} = 2.93 \cdot 10^3 $ \\
	\end{tabular}
	\caption{Terms with parameters in the potential energy of alanine-dipeptide.}
	\label{tab:alaninedipeptide}
\end{table}

For the remainder of this section, we only use the variable $\psi$ as a reaction coordinate, i.e., 
\[
\xi(x) = \psi,
\]
 since $\psi$ is the slowest degree of freedom in the system. 
 In all numerical experiments in this section, we base the approximate macroscopic distribution $\bar{\mu}_0$ on the exact free energy function of $\psi$ and the macroscopic proposal moves are generated by the Euler-Maruyama discretization of the overdamped Langevin dynamics based on this free energy function. We will also fix the inverse temperature at $\beta = 1/100$.

\paragraph{Overview of this section}
We perform three numerical experiments in this section. First, we numerically determine a near-optimal value of $\lambda$ for alanine-dipeptide in Section~\ref{subsubsec:alanine_optimal_lambda} since this analysis should be done on any molecular test case. Then, once we have obtained this near-optimal value, we visually inspect the efficiency of mM-MCMC over the microscopic MALA method by investigating the histogram fits on the distributions of $\psi$ and $\phi$ by both methods. Finally, we compute the total efficiency gain of mM-MCMC on the estimated mean and variance of $\psi$ and $\phi$ and also show some statistics of both numerical sampling methods.

\subsubsection{Optimal value of $\lambda$} \label{subsubsec:alanine_optimal_lambda}
\paragraph{Experimental setup}
Before applying the mM-MCMC method with indirect reconstruction to the alanine-dipeptide molecule, we must first find a near-optimal value for the strength of the biasing potential, $\lambda$. Once we have such this optimal value, we can visually assess the performance of mM-MCMC and compute the efficiency gain over the microscopic MALA method. This analysis should be performed on any molecular system where one wants to use the mM-MCMC method with indirect reconstruction and that is why we start with this analysis. 

\begin{figure}
	\centering
	\includegraphics[width=0.6\linewidth]{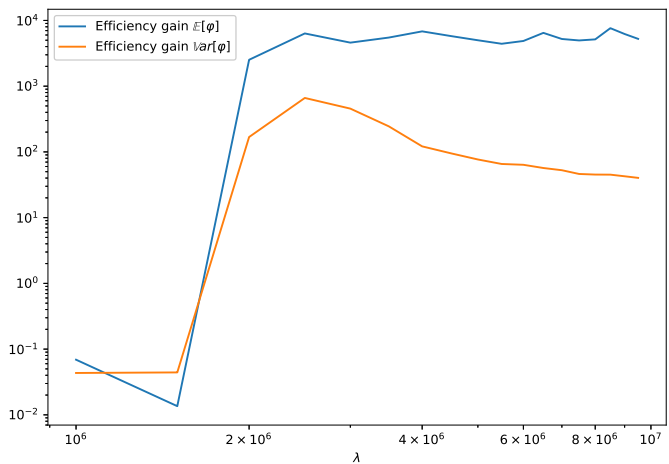}
	\caption{The efficiency gain of nM-MCMC with indirect reconstruction over the microscopic MALA method on the estimated mean (blue) and the estimated variance (orange) of $\psi$. Whenever $\lambda \geq 2 \ 10^6$, the efficiency gain on the estimate mean of $\phi$ is constant, but there is a clear maximal efficiency gain on the estimated variance of $\psi$ at $\lambda = 2.5 \cdot 10^6$.}
	\label{fig:alanine_lambda}
\end{figure}

In this experiment, we define $18$ equidistant values of $\lambda$ between $10^6$ and $10^7$, and for each value, we compute the efficiency gain of mM-MCMC over the microscopic MALA method for $N=10^6$ sampling steps. For the mM-MCMC method, the macroscopic proposals are based on the Euler-Maruyama discretization of the overdamped Langevin dynamics of $\psi$ with time steps $\Delta t = 0.001$, and we take $K=8$ biased steps with time step $0.5 \lambda^{-1}$. The time step for the microscopic MALA method is fixed at $10^{-7}$, while the temperature parameter is $\beta = 1/100$. Figure~\ref{fig:alanine_lambda} displays the efficiency gain of mM-MCMC on the estimated mean and variance of reaction coordinate $\psi$ as a function of $\lambda$.

\paragraph{Numerical results}
Based on Figure~\ref{fig:alanine_lambda}, we see that the efficiency gain on the estimated mean of $\psi$ is constant when $\lambda \geq 2 \ 10^6$. However, when $\lambda$ is smaller than $2 \ 10^6$, the efficiency gain decreases fast beneath $1$. This result is consistent with our claim from Section~\ref{subsec:optimalparameters} that the optimal value of $\lambda$ is slightly larger than the stiffest mode of the system, which is $k_{CC}$ is this case in Table~\ref{tab:alaninedipeptide}. Further, we see that there is a clear optimal value of $\lambda = 2.5 \ 10^6$ that gives the highest efficiency gain on the estimated variance of $\psi$. This results is in contrast to the three-atom molecule where there is a range of optimal values of $\lambda$ on the estimated variance of the reaction coordinate (see Section~\ref{subsubsec:threeatomparam}).
Therefore, we will use this optimal value $\lambda  =2.5 \ 10^{6}$ in the subsequent numerical experiments on the alanine-dipeptide molecule.

\subsubsection{Visual inspection of mM-MCMC} \label{subsubsec:visual_alanine}
\begin{figure}
\centering
\begin{subfigure}[b]{0.5\textwidth}
	\centering
	\includegraphics[width=0.7\linewidth]{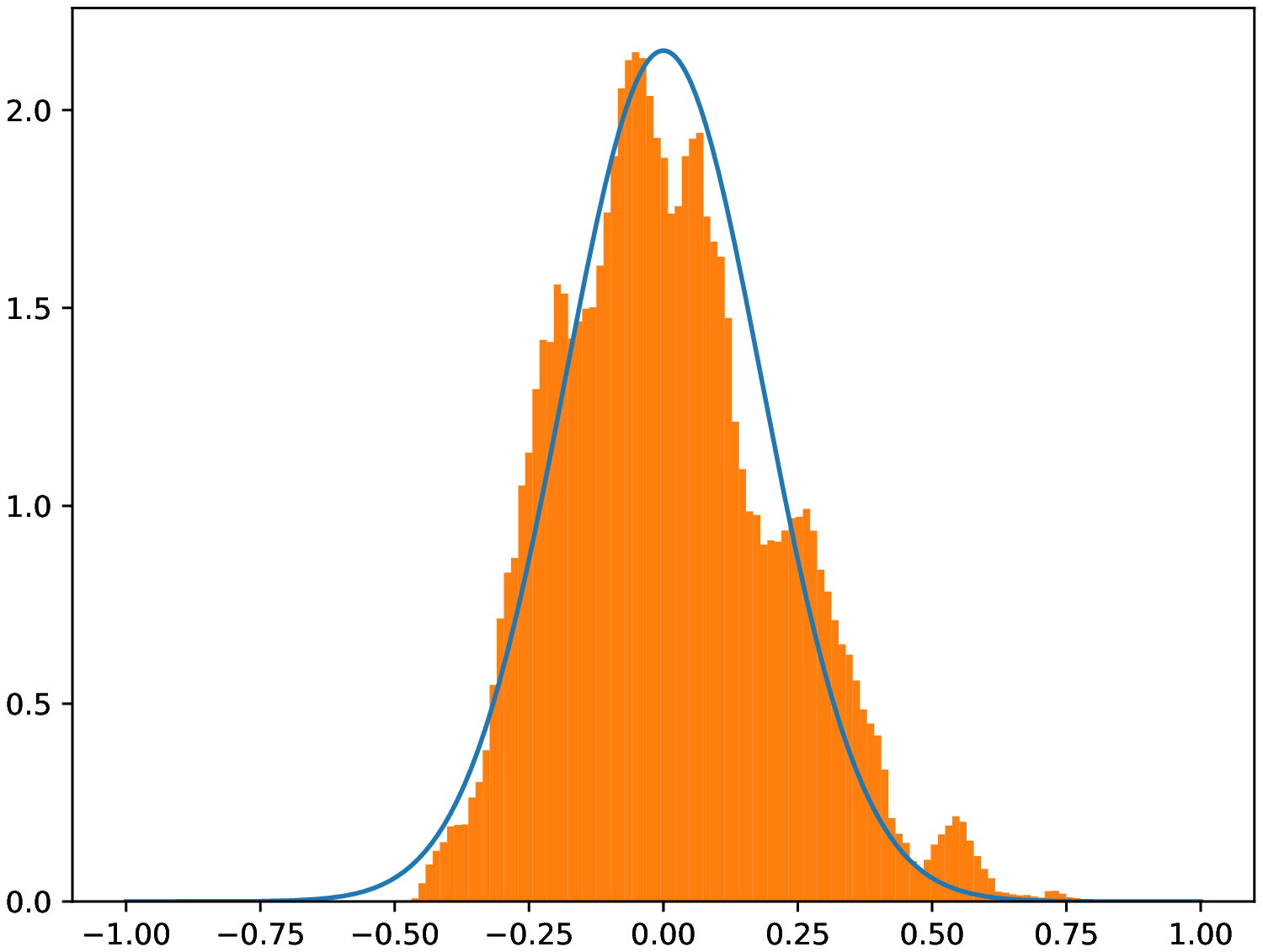}
\end{subfigure}%
\begin{subfigure}[b]{0.5\textwidth}
	\centering
	\includegraphics[width=0.7\linewidth]{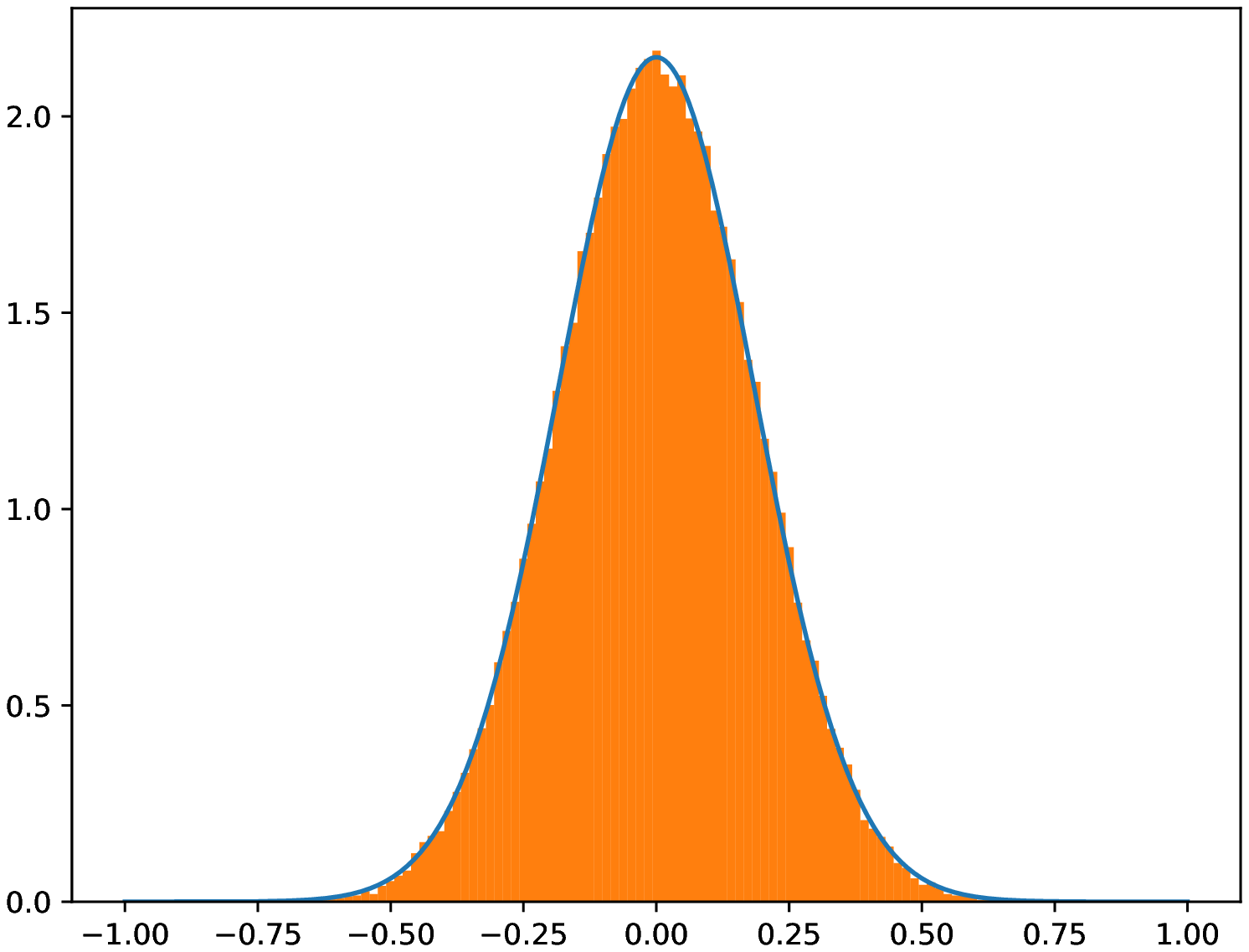}
\end{subfigure}
\begin{subfigure}[b]{0.5\textwidth}
	\centering
	\includegraphics[width=0.7\linewidth]{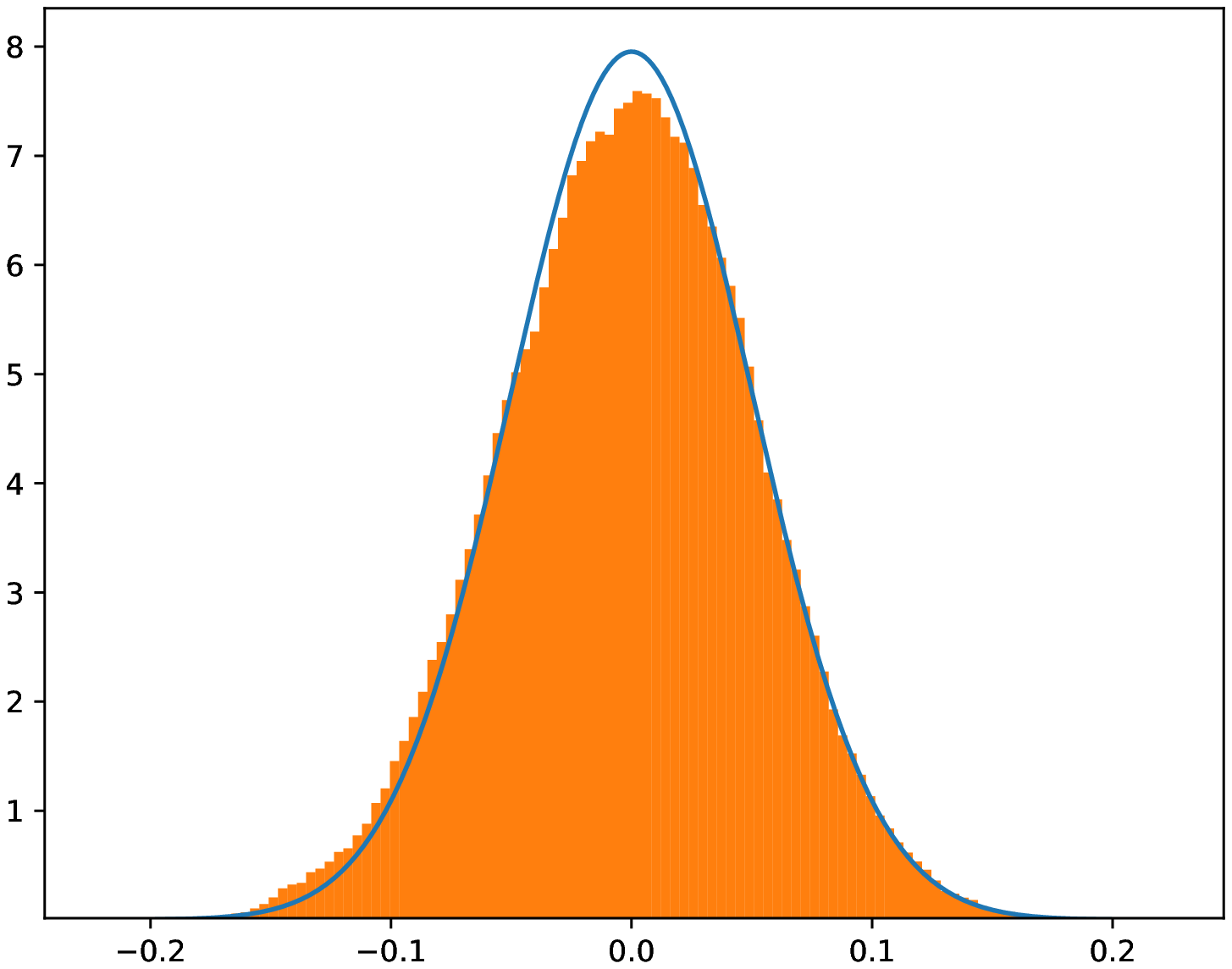}
\end{subfigure}%
\begin{subfigure}[b]{0.5\textwidth}
	\centering
	\includegraphics[width=0.7\linewidth]{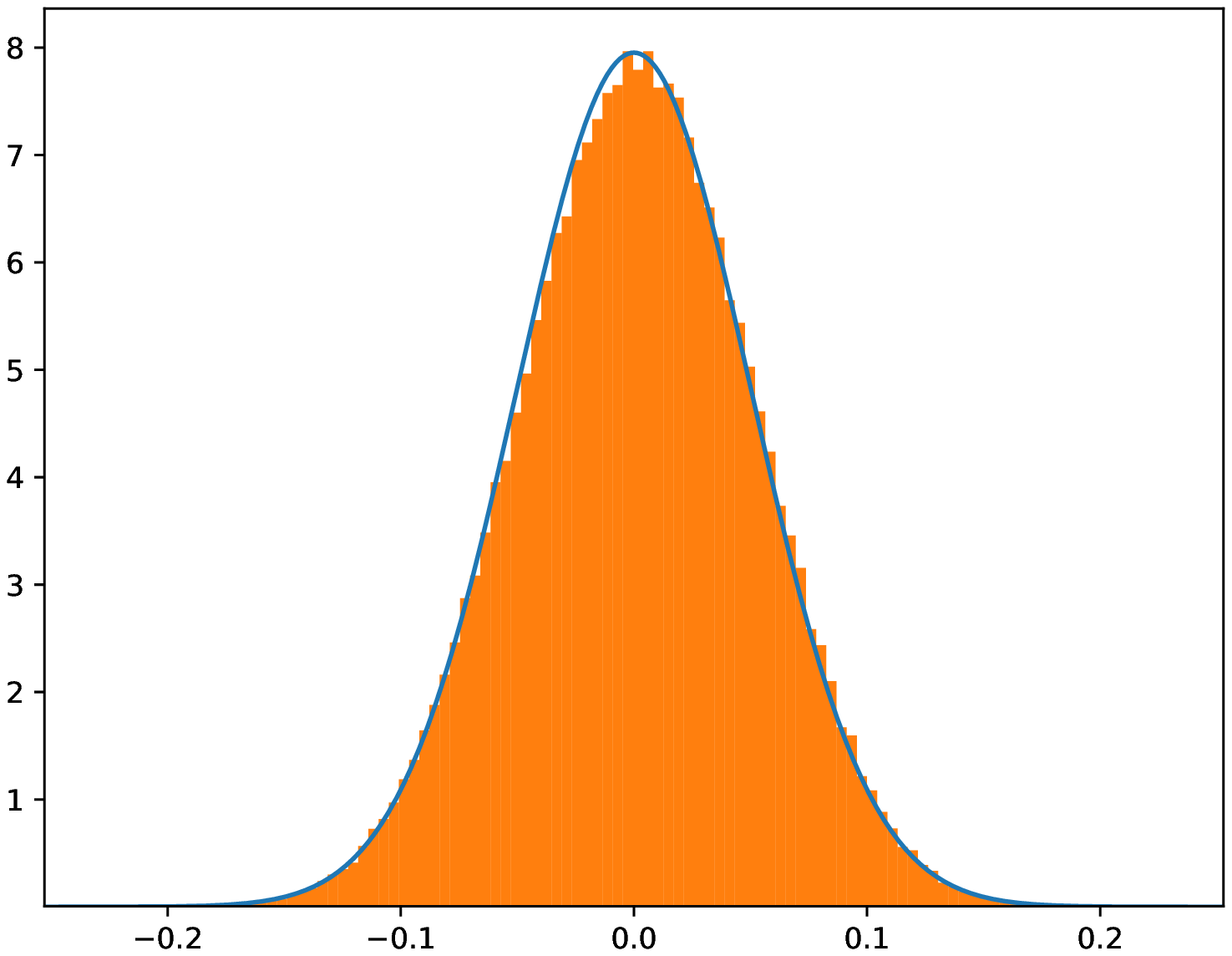}
\end{subfigure}
\caption{Histograms of $\psi$ (top) and $\phi$ (bottom), computed using the microscopic MALA method (left) and mM-MCMC (right). MCMC yields a bad fit on $\psi$ due to the large time-scale separation and the fit on $\phi$ is better. mM-MCMC gives a good fit on both histograms, but we expect the gain on $\psi$ to be larger than the gain on $\psi$.}
\label{fig:hist_alanine}
\end{figure}

\paragraph{Experimental setup}
Now that we have a good estimate for the value of $\lambda$ that yields the optimal efficiency gain, we can visually compare the mM-MCMC method with the microscopic MALA method. Particularly, we consider the histogram fit of the microscopic MALA and the mM-MCMC methods on the marginal distribution of the slow torsion angles $\psi$ and $\phi$. The numerical parameters are $K=8, \ \beta = 1/100, \ \lambda =2.5 \cdot 10^6, \Delta t=0.001$ and for biased simulation time step we take $0.5/\lambda$. The time step of the microscopic MALA method is $10^{-7}$. The histograms for $\phi$ and $\psi$ are shown on Figure~\ref{fig:hist_alanine}.

\paragraph{Numerical results}
Clearly, the histogram fit obtained by the microscopic MALA sampler on the distribution of $\psi$ is inaccurate due the large time-scale separation. The fit on the distribution of $\phi$ is more accurate, since the time-scale separation between $\phi$ and the fast dynamics is smaller than with $\psi$, see Table~\ref{tab:alaninedipeptide}. On the other hand, the mM-MCMC method with indirect reconstruction is able to sample the distribution of $\psi$ more accurately due to the large time steps at the macroscopic level. Additionally, mM-MCMC also samples the distribution of $\phi$ as accurately as the microscopic MALA method because the indirect reconstruction step reconstructs a microscopic sample close to the exact time-invariant direct reconstruction distribution~\eqref{eq:nu} on the sub-manifold of a given value of $\psi$.

\subsubsection{Efficiency gain of mM-MCMC on $\psi$ and $\phi$} \label{subsubsec:efficiency_alanine}\
\paragraph{Experimental setup}
Finally, we measure the efficiency gain the mM-MCMC algorithm with indirect reconstruction over the microscopic MALA method on both the estimated mean and variance of the torsion angles $\psi$ and $\phi$. For a statistically relevant comparison, we average these quantities over 100 independent runs, and we also record the runtime and the variance on the quantity of interest. The numerical results are summarised in Table~\ref{tab:gain_alanine-dipeptide}.

\begin{table}
	\centering
	\begin{tabular}{c|c|c|c|c|c}
		$F$ & Variable & Method & Variance  & Runtime $ \times 10^3$ seconds & Gain \\
		\hline
		\multirow{4}{*}{$\mathbb{E}[\cdot]$}&\multirow{2}{*}{$\psi$}&MCMC &  $1.43 \cdot 10^{-3}$  & $17.5$ & \multirow{2}{*}{$6339.5$} \\
		&&mM-MCMC & $1.25 \cdot 10^{-7}$  & $31.6$&\\
				&\multirow{2}{*}{$\phi$}&MCMC & $7.43 \cdot 10^{-6}$  &  $17.5$& \multirow{2}{*}{$4.22$} \\
		& & mM-MCMC & $9.75 \cdot 10^{-7}$  & $31.6$& \\
		
		\hline
		\multirow{4}{*}{$\text{Var}[\cdot]$}&\multirow{2}{*}{$\psi$}&MCMC &  $6.46 \cdot 10^{-5}$  & $17.5$& \multirow{2}{*}{$661.77$} \\
		&&mM-MCMC & $5.42 \cdot 10^{-8}$  & $31.6$&\\
		&\multirow{2}{*}{$\phi$}&MCMC & $1.31 \cdot 10^{-8}$  &  $17.5$& \multirow{2}{*}{$3.14$} \\
		& & mM-MCMC & $2.31 \cdot 10^{-9}$  & $31.6$& \\
	\end{tabular}
	\caption{Experimental results of the microscopic MALA and mM-MCMC methods for alanine-dipeptide.}
	\label{tab:gain_alanine-dipeptide}
\end{table}

\paragraph{Numerical results}
First of all, there is a large efficiency gain on the estimated mean of $\psi$, as we anticipated due to the slow nature of $\psi$. Additionally, there also is a large gain on the estimated variance of $\psi$, although smaller than the gain on the mean of $\psi$. We do not have an expression for the maximal efficiency gain of mM-MCMC for a given variable, so we cannot compare with the theoretical expressions, but an efficiency gain factor of  $6339.5$ is significant in practice. Finally, there is also a small efficiency gain on the estimated mean and variance on $\phi$. These efficiency gains are smaller than the corresponding efficiency gains on the estimated mean and variance of $\psi$, because $\psi$ is a slower variable than $\phi$, as we noted on Figure~\ref{fig:hist_alanine}. However, as we also mentioned in the previous section, the reason that we obtain an efficiency gain on the estimates of $\phi$ is because the indirect reconstruction distribution lies close to the exact time-invariant reconstruction distribution~\eqref{eq:nu}.

\section{Conclusion and outlook} \label{sec:conclusion}
We presented a micro-macro Markov chain method with indirect reconstruction to sample time-invariant Gibbs distributions from molecular dynamics where there is a time-scale separation between a low-dimensional (macroscopic) reaction coordinate and the high-dimensional (microscopic) degrees of freedom. This method is based on the mM-MCMC method \emph{direct} reconstruction that we worked out in the companion paper~\cite{vandecasteele2020direct}. The mM-MCMC method with direct reconstruction is not always applicable however, since we need to reconstruct a microscopic sample defined on some (possibly) highly non-linear sub-manifold of all microscopic samples that have a fixed value for the macroscopic reaction coordinate. With the indirect reconstruction scheme, we significantly extended the micro-macro Markov chain Monte Carlo (mM-MCMC) method to any reaction coordinate function. The indirect reconstruction step is based on a stochastic process that is strongly biased to a fixed reaction coordinate value. After a few steps of this biased process, the microscopic sample that is obtained has a reaction coordinate value close to the value sampled at the macroscopic level. However, since the indirect reconstruction scheme results in a microscopic sample that does not directly onto the given sub-manifold, we extended the state space such that it includes the dimensions of the reaction coordinate values. The indirect reconstruction scheme also allows us to pre-compute the free energy and the effective dynamics coefficients in a general manner. These computations are consistent with the exact free energy and effective dynamics when the strength of the biased potential increases to infinity. We also investigated the efficiency gain of mM-MCMC with indirect reconstruction on two molecular test cases: a simple three-atom molecule and alanine-dipeptide. On both examples, there is a clear efficiency gain on the order of the time-scale separation.

The indirect reconstruction scheme opens up a series of directions for further research. We mention two main directions here. First, the mM-MCMC method both with direct and indirect reconstruction rely on an available approximation to the free energy function. The free energy is expensive to compute and it would be beneficial if we could evaluate the free energy difference on the fly. Such a scheme would be especially useful in a multilevel scheme where we only use the mM-MCMC scheme to correct for errors of some macroscopic sampler. Second, it would be of interest to investigate the performance the mM-MCMC method with indirect reconstruction on a range of realistic, high-dimensional molecular problems and study whether one can maintain the significant efficiency gain in practical applications.

\bibliographystyle{plain}
\bibliography{refs}

\end{document}